\RequirePackage{fix-cm}
\documentclass{svjour3}                     
\smartqed  
\usepackage{graphicx}
\usepackage{amsmath}
\usepackage{amssymb}
\usepackage{amsfonts}
\usepackage{amsbsy}
\usepackage{color}
\usepackage{mathptmx}      
\usepackage{latexsym}
\begin{document}

\title{A conservative implicit multirate method for hyperbolic problems
}
\subtitle{}


\author{Ludovica Delpopolo Carciopolo \and
		Luca Bonaventura \and Anna Scotti \and Luca Formaggia 
}


\institute{L. Delpopolo Carciopolo \and
		   L. Bonaventura \and A. Scotti \and L. Formaggia
           \\
          MOX, Dipartimento di Matematica,\\
          Politecnico di Milano, Via Bonardi 9, 20133 Milano, Italy\\
           \and
           L. Delpopolo Carciopolo\\
           \email{ludovica.delpopolo@polimi.it}   \\
           \and
		   L. Bonaventura\\
		   \email{luca.bonaventura@polimi.it}   \\
		   \and
           A. Scotti\\
           \email{anna.scotti@polimi.it}   \\
           \and
           L. Formaggia\\
           \email{luca.formaggia@polimi.it}   \\
}

\date{\today}

\maketitle

\begin{abstract}
This work focuses on the development of a self adjusting multirate strategy based on  an implicit time discretization for the numerical solution of  hyperbolic equations, that could benefit from different time steps in different areas of the spatial domain. We propose a novel   mass conservative multirate approach, that can be generalized to various implicit time discretization methods. It is based on flux partitioning, so that flux exchanges between a cell and its neighbors are balanced.
 A number of numerical experiments on both non-linear scalar problems and systems of hyperbolic
 equations have been carried out to test the efficiency and accuracy of the proposed approach.
\keywords{Multirate schemes, Conservation laws, Conservative formulation.}
\end{abstract}

\section{Introduction}
 \label{intro} \indent
Conservation laws model a large variety of phenomena in the geosciences, such as shallow water flow, multiphase groundwater flows and advection and dispersion of contaminants.
The time discretization of hyperbolic problems is often subject to restrictions on the time step. Explicit time integration schemes are only stable if the time step amplitude fulfils the well known CFL condition \cite{leveque:1992}, an upper bound dictated by the space discretization parameter and the wave speed. Thus, a small mesh size or a high wave speed in a small part of the domain imposes a strict limitation on the time step everywhere. To overcome this problem, it is possible to use  implicit, unconditionally stable methods, which allow larger time steps, but require the solution of a possibly non linear system at each time step. Moreover, all high order implicit scheme are not unconditionally monotone, so that a different condition on the size of time step is required to ensure monotonicity. Finally, in the case of systems representing phenomena evolving on multiple time scales,
implicit schemes allow to approximate correctly the slower components of the solution only at the price
of a significant loss of accuracy on the faster ones.
For these reasons, we  investigate in this paper the benefits of a multirate approach for these  problems.
 
Multirate methods were originally proposed in \cite{rice:1960} in the context of systems of ordinary differential equations. Many studies  have been then devoted to the improvement of these methods, see e.g. \cite{andrus:1979}, \cite{gear:1984}.   The main idea of multirate methods is to integrate each component of the system using a different time step. Slow components, i.e. components with longer characteristic time scales, are integrated with larger time steps, while smaller time steps are used only for fast components. Thus,   multirate methods   can avoid a significant amount of the computations that are necessary in the single rate approaches, if the faster components that require a small time step are confined in a small part of the domain (possibly evolving in time). In other words, in the multirate approach the most appropriate time resolution is employed for each variable of the system.
 In   earlier multirate methods,
the system was partitioned  \textit{a priori}, based on the knowledge of the specific problem to be solved.  
 A self adjusting, recursive time stepping strategy has been then proposed in  \cite{savcenco:2007}.
In this more recent approach, a tentative  global  step is first taken for all components, using a robust, unconditionally stable method. The time step is then reduced
only  for those components for which a suitable local error estimator is greater than the specified tolerance.
In this way, automatic detection of fast components is achieved.

In   \cite{Constantinescu:2007} and  \cite{fok:2015}
the authors propose multirate Runge-Kutta methods that preserve the stability properties
 of the single rate approach. 
 We will base our work on the strategy proposed in
\cite{savcenco:2007} for the $\theta$-method and extended in \cite{bonaventura:2018} to the TR-BDF2 method  as fundamental single rate solver. The TR-BDF2 method has been originally introduced in \cite{bank:1985} and more thoroughly analyzed in  \cite{hosea:1996}. It is a second order, one step,  L-stable implicit method endowed with a number  of interesting properties, as discussed in  \cite{hosea:1996}. As in \cite{bonaventura:2018}, in our approach 
the choice of the time step size at each step is based on the technique proposed in~\cite{fok:2015}.

While multirate methods have been mostly applied to general systems of ODEs, in this work we will 
focus exclusively on   systems  that arise from the space discretization of conservation laws. Unlike previous attempts, we propose a component partitioning strategy which is based on the the numerical fluxes, in order to preserve the mass conservation properties of the single rate method.  This approach is inspired by the flux partitioning strategy proposed in \cite{ketcheson:2013} and already successfully employed in \cite{bonaventura:2017} to derive monotonic methods for space discretized conservation laws.

This paper is structured as follows. 
In Sect. \ref{selfa}, the multirate approach of \cite{bonaventura:2018} is briefly reviewed.
 In Sect. \ref{conservation} we describe in detail  the conservative  algorithm and we present
 a brief analysis on the consistency of the method. In Sect. \ref{numerical} numerical results for nonlinear conservation laws are presented. Conclusions are drawn in the final section.

  \section{A self adjusting   multirate approach}
 \label{selfa} 
 \indent
 In this section, the self adjusting multirate approach proposed in \cite{bonaventura:2018}
 is outlined, as applied to   the solution of the Cauchy initial value problem
\begin{equation}\label{eq:cauchy}
y^{\prime}(t)= f(t, y(t)),\quad t\in (0,T],\quad  y(0)= y_0\in\mathbb{R}^m.
\end{equation}
We   consider time discretizations associated to discrete time levels
$t_n,\  n=0,\dots,N $ such that $\Delta t_n=t_{n+1}-t_n$ and
we will denote by
$u^n $ the numerical approximation of $y(t_n).$ 
 We will also denote by $ u^{n+1}={\cal S} (u^{n},\Delta t_n) $ the
implicitly defined operator ${\cal S} :\mathbb{R}^m\rightarrow \mathbb{R}^m $
 whose application is equivalent to the computation of  one step of size $\Delta t_n$  of
a given single step method. While here only implicit methods will be considered,
the whole framework can be extended to explicit and IMEX methods.
Notice that, if $ P$ is the projector onto a linear subspace ${\cal V} \subset {\mathbb{R}}^m $ with dimension $p < m,$
the operator ${\cal S}^{\cal V} :\mathbb{R}^p\times \mathbb{R}^{m-p}\rightarrow \mathbb{R}^p $
that represents the solution of the subsystem obtained freezing the components
of the unknown belonging to  ${\cal V}^{\perp}  $  to the value $z \in \mathbb{R}^{m-p}$ can be defined by
$ y =   {\cal S}^{\cal V} ( x, z, \Delta t_n )=  P{\cal S} ( x\oplus z, \Delta t_n) .$
Furthermore, we will denote by $ Q(u^{n+1},u^{n},\zeta) $ 
the interpolation operator that provides an approximation of the numerical solution at
intermediate time levels $t_n+\zeta,$ where $\zeta \in[0,\Delta t_n].$ Linear interpolation is often employed, but, for
 multistage methods, knowledge of the intermediate stages  also allows the application of more accurate interpolation procedures  without 
substantially increasing the computational cost.

In a multirate approach,   system \eqref{eq:cauchy} is partitioned into a sub-system of so called \emph{active components} with a faster time scale and  into the complementary sub-system of the \emph{latent components}, which are associated to slower phenomena.
In this context, the basic idea of a self-adjusting strategy is to use a tentative global time step to identify the set of the active components, which have to be recomputed with a smaller time step to maintain the desired accuracy and stability. In particular, the self-adjusting  multirate algorithm introduced in \cite{bonaventura:2018}   is a generalization of that proposed in \cite{savcenco:2007} and can be described as follows.
\begin{itemize}
	\item[1)]  Perform  a tentative global (or macro) time step of size $\Delta t_n $ with the standard single rate method and compute
	 $\hat u^{n+1}= {\cal S} (u^{n},\Delta t_n)$.
	 \item[2)]  Apply the error estimator 
	 to  partition the state space into  active and  latent variables.
	  The projection onto the subspace ${\cal V}_{0}$
	of the active variables is denoted by $P_n^{(0)},$  
	while the projection onto the complementary
	subspace will be denoted by $\bar{ P}_n^{(0)}.$
	Define $    \bar{P}_n^{(0)} u^{n+1} =    \bar{P}_n^{(0)}\hat { u}^{n+1} $
	as well as $ {u}^{n,0}= {u}^{n} $ and
	  $t_{n,0}=t_{n}. $ 
	\item[3)]  For $k\geq 1,$
	choose a local (or micro) time step  $\Delta t_n^{(k)} $ for the active variables, based on the value of the error estimator.  Set  $t_{n,k}=\min\{t_{n,k-1}+\Delta t_n^{(k)},  t_{n+1}\}. $ 
	\begin{itemize}
		\item[3.1)]      Update the latent variables by interpolation
	$$	
	 \bar{ P}_n^{(k)}{ u}^{n,k}=
	 Q  ( \bar{ P}_n^{(k-1)}{u}^{n+1},  \bar{P}_n^{(k-1)} { u}^{n,k}, \Delta t_n^{(k)} ).
	$$
	\item[3.2)]  Update the active variables  by computing
		$$
	 P^{(k)}_n{u}^{n,k}=
	  {\cal S}^{{\cal V}_{k-1}} (P^{(k-1)}_n{u}^{n,k},  \bar{P}_n^{(k-1)}{u}^{n,k-1}, \Delta t_n^{(k)}).
	$$
	\item[3.3)] Compute the error estimator  for the active variables only and   partition again  $ {\cal V}_{k-1} $
	into latent and active variables. Denote by ${\cal V}_{k}\subset {\cal V}_{k-1} $ the new
	subspace of active variables and  by  $ P_n^{(k)} $  the corresponding projection.	
	\item[3.4)] Repeat  3.1) - 3.3) until $t_{n,k}=t_{n+1}.$
 
		\end{itemize}
	\end{itemize}

 A stability analysis of the above described approach has been
proposed in \cite{bonaventura:2018} in the case of a linear system with a simplified refinement strategy.
The effectiveness of the above procedure depends in a crucial way on the accuracy and stability of the basic ODE 
solver ${\cal S}, $ as well as on the time step
refinement and partitioning criterion.  In  \cite{bonaventura:2018}, the embedded error estimator of the TR-BDF2 method was used for the error estimator and
 the error control strategy  proposed in  \cite{fok:2015} was extended to  employ  a combination of absolute and relative error tolerances.  It is important to remark that the previously defined approach,
 when applied to ODE systems stemming from the space discretization of conservation laws like \eqref{eq:semidiscr},
  does not  guarantee mass conservation for the numerical solution, since some of the fluxes
  are recomputed during refinement only for one of the two adjacent variables. For this reason,
  in section \ref{conservation} we propose a conservative version of the method.
 
 \section{The conservative implicit multirate approach for hyperbolic conservation laws}
\label{conservation}
The aim of this section is to introduce a mass conservative, implicit multirate scheme to integrate in time non linear conservation laws of the form
\begin{equation*}
\frac{\partial u}{\partial t} + \frac{\partial f(u)}{\partial x} = 0 \qquad x \in\mathbb{R}, \quad t>0,
\end{equation*}
with given initial datum $u(x,0) = u_{0}(x)$ for $ x \in \mathbb{R}$.
For simplicity, in this section, we will treat scalar problems in one-dimension and we pose the differential problem on the whole real line, postponing to a later stage a discussion on how to treat boundary conditions for problem in a bounded domain.
To discretize the equation in space we consider the set of the cells $I_{i} =\left[ {x_{i-\frac{1}{2}}},{x_{i+\frac{1}{2}}} \right] $, for $i\in\mathbb{Z}$,  with ${x_i}$ being the center of cell $I_i$ and $\Delta x_i = {x_{i+\frac{1}{2}}}-{x_{i-\frac{1}{2}}}$ the cell size. 

We denote by $u_i(t)$ the approximation of the average value of $u(x,t)$ in cell $I_i$ after the spatial discretization, i.e
$u_i(t) \simeq \dfrac{1}{\Delta x_i} \int_{x_{i-\frac{1}{2}}}^{x_{i+\frac{1}{2}}} u(x,t) \ dx$ for $t>0$, while 
the initial value at $t=0$ is obtained from the initial data, 
\begin{equation*}
u_i(0) = \dfrac{1}{\Delta x_i} \int_{x_{i-\frac{1}{2}}}^{x_{i+\frac{1}{2}}} u_0(x) \ dx.
\end{equation*}

A conservative finite volume discretization yields the following system of ordinary differential equations
\begin{equation}\label{eq:semidiscr}
 \dfrac{du_i}{dt}(t)= - \frac{1}{\Delta x_i}\left[F_{i+\frac{1}{2}}(t) - F_{i-\frac{1}{2}}(t)\right ],\quad i \in \mathbb{Z},\quad t>0,
\end{equation}
where $F_{i \pm \frac{1}{2}}(t) = F(u_{i\mp p}(t),\cdots,u_i(i),\cdots,u_{i\pm q}(i))$ is the semi-discrete numerical flux at the control volume 
face $x_{i \pm \frac{1}{2}}$ and $x_{i \mp p},\cdots,x_i,\cdots,x_{i \pm q}$ is the stencil of nodes used to evaluate it. For instance, 
in the classical two-point flux approximation $p=0$ and $q=1$. 

Equations in the form \eqref{eq:semidiscr} are the starting point for our multirate approach, which, differently from the scheme outlined in the previous section, employs an error estimator based on the fluxes rather than on the system components to identify active and latent components, with the aim to maintain the mass conservation properties of the basic scheme.

We give here a general overview of the method, postponing to a later section a more detailed description of the algorithm. Given the numerical solution at time $t_n$ and a global time step $\Delta t_{n}^{(0)} = t_{n+1} - t_{n}$,  we aim to a numerical scheme that may be eventually written in the form
\begin{equation}\label{eq:consform}
u_i^{n+1}=u_i^n -\frac{1}{\Delta x_i}\left( H_{i+\frac{1}{2}}-H_{1-\frac{1}{2}}\right)
\end{equation}
where
\[
H_{i \pm\frac{1}{2}}\cong \int_{t_n}^{t_{n+1}} F_{i\pm \frac{1}{2}}\, dt 
\]
is the numerical flux, which typically depends on $F_{i\pm \frac{1}{2}}$ sampled at different times. Note that we are using a 
non-standard definition for the numerical flux, since we are not dividing the time integral by the time step length. Discretizations of the form~\eqref{eq:consform}
are conservative in the sense that, for any set of indices $\mathcal{I},$ the quantity $\sum_{i\in\mathcal{I}}  \Delta_i (u_i^{n+1}-u_i^n)$
depends only on the values of the numerical fluxes at the boundary of the set $\cup_{i\in\mathcal{I}} I_i$.

At each time step, we first compute the approximate solution at $t_{n+1}$ for all components with a tentative time step. The value of the numerical fluxes at all interfaces is checked using an appropriate error estimator. If the flux is rejected on the basis of the error estimator, all components involved in its stencil are added to the set of active components that need to be recomputed with a smaller time step. 
During the re-computation,  the accepted numerical fluxes  are kept constant inside the time slab and interpolation is used to obtain their appropriate value, while the rejected ones are recomputed. 
In this way,  interpolation is  applied directly to the fluxes, rather than to the components, which allows to maintain the structure of the scheme 
in the form~\eqref{eq:consform}, where the $H_{i\pm \frac{1}{2}}$  will consist, at the end of the procedure, of contributions coming from the accepted fluxes.

\subsection{A first example}
\label{example}
For the sake of clarity, we first present the proposed multirate method using the $\theta$-method as implicit time integration scheme, while to discretize in space we adopt a uniform grid with step size $\Delta x$. The purpose is to give an idea of the scheme on a simple example, before presenting the general procedure. We will also assume to employ a two-point flux approximation, which means
$F_{i\pm \frac{1}{2}}=F(u_i,u_{i\pm 1})$ . At the global time level $t_{n}^{(0)}, $ using the time 
step $\Delta t_{n}=\Delta t$, the following expression is obtained in the first tentative calculation:
\begin{equation*}
\hat{u}_{i}^{n+1} = u_{i}^{n} -\frac{\theta}{\Delta x} \left[F_{i+\frac{1}{2}}^{n+1} - F_{i-\frac{1}{2}}^{n+1}\right] -\frac{1-\theta}{\Delta x}\left[F_{i+\frac{1}{2}}^{n} - F_{i-\frac{1}{2}}^{n}\right],
\end{equation*}
where $F^n_{1\pm \frac{1}{2}}$ denotes the numerical flux computed using the value of the approximated components at time $t_n$.
Clearly, with  a simple manipulation the scheme can be rewritten in form~\eqref{eq:consform}. We also  
notice that here  ${\Delta t_n}$ is included in the numerical fluxes, in contrast with other description of the scheme found in the literature.

\begin{figure}
\centering
\includegraphics[scale=0.75]{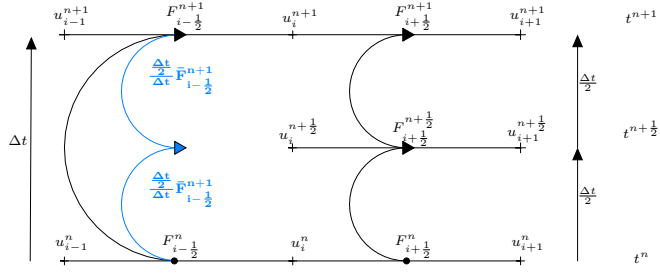}

\caption{An example of flux partitioning that preserves   mass at each global time step.}
\label{fig::tikz_flux}
\end{figure}

If we suppose, as showed in Fig. \ref{fig::tikz_flux}, that at this level the error estimator rejects the flux at the interface point $x_{i+\frac{1}{2}}$, we have to recompute 
the components of the stencil of  $F_{i+\frac{1}{2}}$, i.e. $u_i$ and $u_{i+1}$ will be recomputed  using a smaller time step. Here, for simplicity, we reduce $\Delta t_n^{(0)}$ by a half. If instead  $F_{i-\frac{1}{2}}$ is  accepted, at  the new intermediate time $ t_{n+\frac{1}{2}}=t_n+\Delta t_n^{(1)}=t_n+\frac{1}{2}\Delta t_n^{(0)}$ we have
\begin{equation*}
u_{i}^{n+\frac{1}{2}} = u_{i}^{n} -\frac{\theta}{\Delta x}  \left[F_{i+\frac{1}{2}}^{n+\frac{1}{2}} - \frac{1}{2}{F}_{i-\frac{1}{2}}^{n+1}\right] -\frac{1-\theta}{\Delta x}\left[F_{i+\frac{1}{2}}^{n} - \frac{1}{2}{F}_{i-\frac{1}{2}}^{n}\right].
\end{equation*}
Here, $F^n_{1-\frac{1}{2}}$  and $F^{n+1}_{1-\frac{1}{2}}$ have been kept frozen at the value computed at the larger time step (since $F_{i-\frac{1}{2}}$ has been accepted). They are multiplied by a factor $\frac{1}{2}$ to  account for the time step reduction $\frac{\Delta t_n^{(1)}}{\Delta t_n^{(0)}}$.
As for cell ${i+1}$, if we suppose to accept the numerical flux at the interface point $x_{i+\frac{3}{2}}$, a similar expression is obtained:
\begin{equation*}
u_{i+1}^{n+\frac{1}{2}} = u_{i+1}^{n} -\frac{\theta}{\Delta x}  \left[\dfrac{1}{2}{F}_{i+\frac{3}{2}}^{n+1} - F_{i+\frac{1}{2}}^{n+\frac{1}{2}}\right] -\frac{1-\theta}{\Delta x}\left[ \dfrac{1}{2}{F}_{i+\frac{3}{2}}^{n} -F_{i+\frac{1}{2}}^{n} \right].
\end{equation*}
If the new time step $\Delta t^{(1)}_n$ is such that all fluxes are accepted, we can recompute the solution at time $t_{n+1}$ as
\begin{align*}
u_{i}^{n+1} &= u_{i}^{n+\frac{1}{2}} -\frac{\theta}{\Delta x} \left[F_{i+\frac{1}{2}}^{n+1} - \dfrac{1}{2}{F}_{i-\frac{1}{2}}^{n+1}\right] -\frac{1-\theta}{\Delta x}\left[F_{i+\frac{1}{2}}^{n+\frac{1}{2}} - \dfrac{1}{2}{F}_{i-\frac{1}{2}}^{n}\right],\\
u_{i+1}^{n+1} &= u_{i+1}^{n+\frac{1}{2}} -\frac{\theta}{\Delta x} \left[\dfrac{1}{2}{F}_{i+\frac{3}{2}}^{n+1} - F_{i+\frac{1}{2}}^{n+1}\right] -\frac{1-\theta}{\Delta x} \left[ \dfrac{1}{2}{F}_{i+\frac{3}{2}}^{n} -F_{i+\frac{1}{2}}^{n+\frac{1}{2}}\right].
\end{align*}
For cell $i-1$, if also the flux $F_{i-\frac{3}{2}}$ has been accepted, the solution at time $t_{n+1}$ is simply:
\begin{equation*}
u_{i-1}^{n+1} = u_{i}^{n} -\frac{\theta}{\Delta x} \left[F_{i-\frac{1}{2}}^{n+1} - F_{i-\frac{3}{2}}^{n+1}\right] -\frac{1-\theta}{\Delta x}\left[F_{i-\frac{1}{2}}^{n} - F_{i-\frac{3}{2}}^{n}\right].
\end{equation*}

One can verify that mass conservation at the global step is guaranteed,
 since all fluxes at interface $i+\frac{1}{2}$ and $i-\frac{1}{2}$ cancel each other exactly.
Since the choice of $i$ is arbitrary, this fact holds true for all interfaces.

\subsection{The time refinement and time stepping strategy}
\label{timestep}

We now present the general algorithm to perform numerical integration inside one global step $t_n\to t_{n+1}$. The algorithm is recursive and, inside the global step, we define a new sub-step each time a flux has been rejected at the current sub-step. Moreover, in the general case the refinement ratio can be different from $\frac{1}{2}$. 
We will generically indicate with $\mathcal{A}_C$ and $\mathcal{A}_F$ the set of active components (i.e. those that have to be recomputed) and that of accepted fluxes, respectively. Superscripts may be added to indicate different instances. These sets  always satisfy the   property
\[
\mathcal{A}_C=\lbrace u_i:\ F_{i-\frac{1}{2}}\not \in \mathcal{A}_F\vee F_{i+\frac{1}{2}}\not\in \mathcal{A}_F \rbrace.
\]
We also introduce the vector $\mathcal{T}_F,$ that for each flux in $\mathcal{A}_F$ records the length of  the sub-step at the moment in which the flux has been accepted.
For consistency of notation, we will use subscripts of the form $i\pm \frac{1}{2}$ to indicate fluxes or flux related quantities. We assume that an error estimator for the fluxes is provided and we only consider a two-point flux approximation, although the procedure can be extended to other types of numerical flux constructions.
    
We denote with $S$ the operator that starting from $u^*$ returns the vector 
$u^{\triangle}$ of updated active components within a given sub-step, and also computes the new sets $\mathcal{A}_F$ and $\mathcal{A}_C$, together with the new time step to be used for the refined sub-steps or the next step.

Algorithm $S$ is the building block for the operator $M,$ that is used recursively to compute a single global time step with our multirate method.
It basically takes as input a set of components $u^*$ and a time step, and proceeds recursively across all rejected sub-steps 
to produce the final value at the end of the time step. The parameter $p$ takes track of the level of refinement.
The first time that   the multirate algorithm is applied, $p$ will be equal to $0$, $u^* = u^n$, $u^\triangle = \hat{u}^{n+1}$ and $\Delta t^* = \Delta t_n$.
\begin{center}
\emph{ALGORITHM $M(u^*,\Delta t^* ,p;u^\triangle, \Delta t^\triangle)$
}
\end{center}
\begin{itemize}
\item set $s = 1$;
\item while $t^* + \Delta t^* \leq t^\triangle$ where $t^*$ and $t^\triangle$ are the times where $u^*$ and $u^\triangle$ have been computed, respectively; 
\begin{enumerate}
\item $u^{(s)} = u^*$;
\item if $p=0$ set $\mathcal{A}^{(0)}_C$ equal to the set of all components, $\mathcal{A}^{(0)}_F=\emptyset$
$\mathcal{T}^{(0)}_F=\emptyset$;
\item Call
\begin{multline*}
S(u^*,\mathcal{A}^{(p)}_C,\mathcal{A}^{(p)}_F,\mathcal{T}^{(p)}_F,\Delta t^*;
u^\triangle,\mathcal{A}^{(p+1)}_C,\mathcal{A}^{(p+1)}_F,\mathcal{T}^{(p+1)}_F,\Delta t^\triangle)
\end{multline*}
\item if $\mathcal{A}^{(p+1)}_C\ne \emptyset$
\begin{itemize}
\item $M(u^{(s)},\Delta t^\triangle,p+1;u^{(s+1)})$;
\end{itemize}
\item otherwise
\begin{itemize}
\item set $u^* = u^\triangle$ and so $t^* = t^\triangle$;
\item set $s = s + 1$;
\end{itemize}
\end{enumerate}

\end{itemize}
The index $p$ indicates the level of refinement, while the index $s$ is the sub-step taken at each level of refinement.
Note that the set of fluxes marked as accepted at the given level are kept as such on all sub-steps associated to that level.  This is the key for mass conservation, as explained later.

The operator $S$ is defined by the following algorithm

\begin{figure*}
\centering
\includegraphics[scale=0.4]{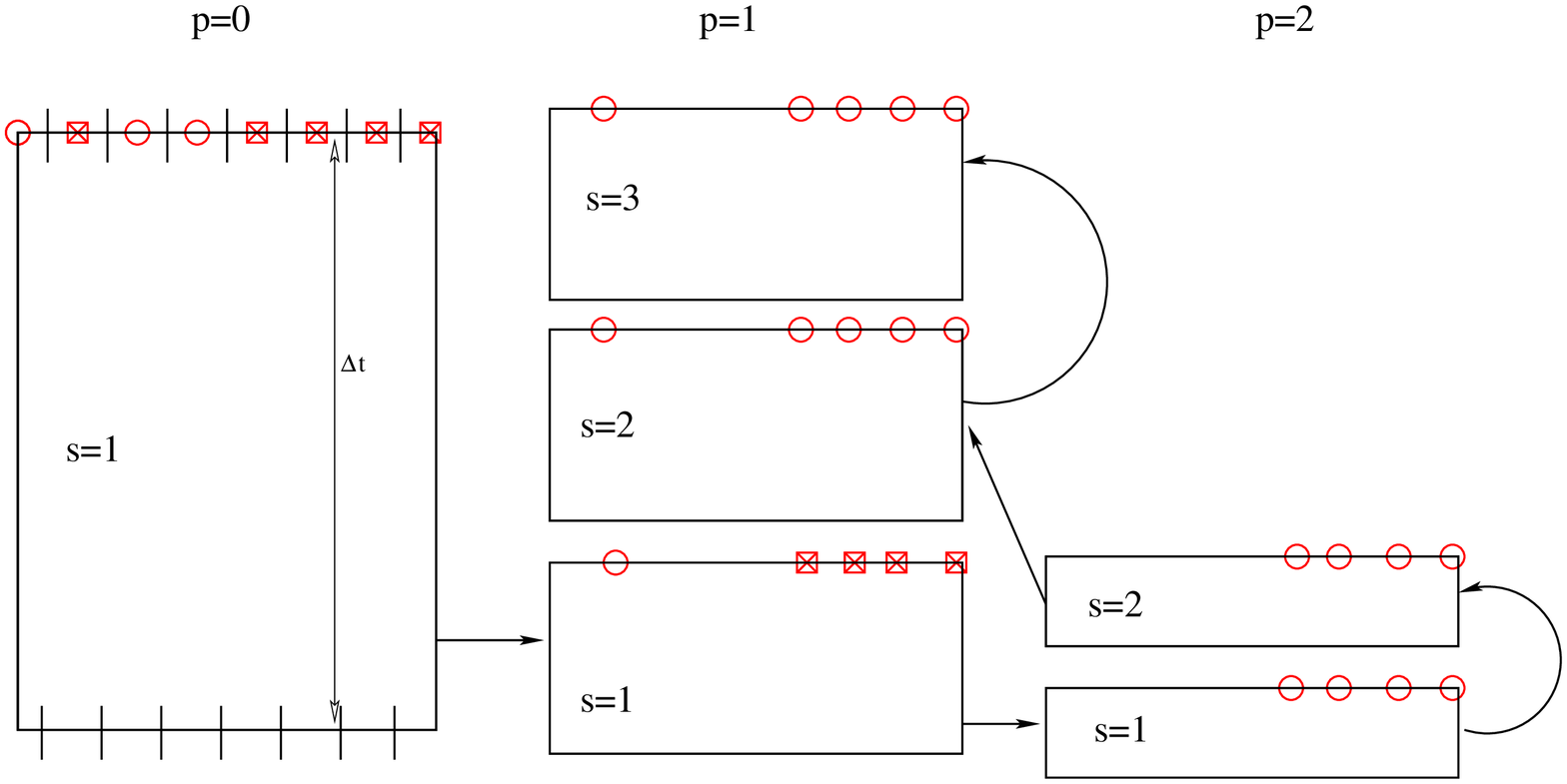}
\caption{Example of the time stepping multirate method.}
\label{fig::tikz_diagram}
\end{figure*}

\begin{center}
\emph{ALGORITHM $S(u^*,\mathcal{A}^*_C,\mathcal{A}^*_F,\mathcal{T}^*_F,\Delta t^*;u^{\triangle},\mathcal{A}^\triangle_C,\mathcal{A}^\triangle_F,\mathcal{T}^\triangle_F,\Delta t^\triangle)$}
\end{center}
\begin{enumerate}
\item Compute $u^{\triangle}$ for all components in $\mathcal{A}^*_C$ starting from  $u^*$ and using the chosen time-advancing scheme with time step
$\Delta t^*$. The fluxes necessary for this computation are given by 
$\mathcal{F}=\lbrace F_{i\pm \frac{1}{2}}:\, u_i^*\in\mathcal{A}^*_C \rbrace$. 
Those contained in $\mathcal{A}^*_F$ will not be recomputed but used 
directly, scaled by the factor $\Delta t^*/\Delta t_{i+\frac{1}{2}}$, where $\Delta t_{i+\frac{1}{2}}$ indicates the corresponding element of $\mathcal{T}_F^*$;
\item Estimate the  error $\epsilon_{i+\frac{1}{2}}$ on all recomputed fluxes, i.e. the fluxes in $\mathcal{F}\setminus\mathcal{A}_F$, to identify the set $\mathcal{R}_F$ of rejected fluxes at this level, $\mathcal{R}_F=\lbrace F_{i+\frac{1}{2}}:\ \epsilon_{i+\frac{1}{2}}> \text{tol} \rbrace$, where $\text{tol}$ is a given tolerance parameter.
\item If $\mathcal{R}_F\ne \emptyset$ compute 
\begin{itemize}
\item The set of active component for the next substep
\[
\mathcal{A}^\triangle_C=\lbrace u_i:\ F_{i-\frac{1}{2}}\in \mathcal{R}_F\vee F_{i+\frac{1}{2}} \in \mathcal{R}_F \rbrace;
\]
\item  The time step for the active components to be recomputed at the next sub-step. We adopt this extension of the formula originally proposed  \cite{fok:2015} and already adapted in  \cite{bonaventura:2018}
$$  \Delta t^{new}= \nu \min_{F_{i+\frac{1}{2}} \in \mathcal{R}_F} \left( \frac{\tau_r |F_{i+\frac{1}{2}}|+ \tau_a}{\epsilon_{i+\frac{1}{2}}} \right)^{\frac{1}{r+1}},$$
where $\tau_r$ and $\tau_a$ are a relative and absolute tolerance, respectively, $r$ is the order of convergence of the chosen time advancing method and $\nu$ an user defined 
parameter taking values in $[0,1].$ As customary in  adaptive time integration approaches, see e.g.
\cite{prince:1981}, these  parameters are employed to tune the adaptation criterion and to impose  a more
conservative choice of the time step if necessary.
\item Set $\Delta t^\triangle$ as the nearest fraction of $\Delta t^*$ smaller than $\Delta t^{new}$;
\end{itemize}
\item Otherwise, set $\mathcal{A}^\triangle_C=\emptyset$ and $\Delta t^\triangle = \Delta t^*$;
\item Return in $\mathcal{A}^\triangle_F$ the set of accepted fluxes for the next level, by setting $\mathcal{A}^\triangle_F=\mathcal{F}\setminus\mathcal{R}_F$, as well as the corresponding $\mathcal{T}^\triangle_F$ for the next level: for the fluxes in $\mathcal{F}\setminus\mathcal{R}_F$ that had already been accepted we just copy the previous value, for the newly accepted fluxes we set it equal to $\Delta t^*$.
\end{enumerate}
We mention that the algorithm keeps track also of the time instants the fluxes have to be computed, for the sake of simplicity we have omitted to indicate it explicitly. In Fig. \ref{fig::tikz_diagram} we draw an example of what it is obtained combining the two algorithms, the circles indicate the latent components inside the sub-step, instead the crosses indicate the active components that have to be recomputed in the next sub-refinement.

\subsection{Mass conservation}
Given a time step $\Delta t_{n} = t_{n+1}-t_n$, the values of the numerical approximation can be written as:
\begin{equation}
u_i^{n+1} = u_i^n  - \frac{1}{\Delta x}\left[ H^+_{i+\frac{1}{2}} -H^-_{i-\frac{1}{2}} \right].
\end{equation}
The fluxes to the right and to the left of each cell are marked with the superscripts $+$ and $-$, respectively, because they have an apparent dependence on the considered cell. 
The aim of this section, however, is to prove that, given two cells, $i$ and $i-1$ for example, the flux at the common interface has the same value $H^-_{i-\frac{1}{2}} = H^+_{i-\frac{1}{2}},$ in spite of being computed by an apparently different procedure.
The fluxes, using the $\theta$-method as time-advancing scheme, can be written in the following way:
\begin{align*}
H^+_{i-\frac{1}{2}}  &= 
\sum_{p,s: \{ u_i^{s,p} \in \mathcal{A}_C^{p,s*} \wedge u_i^{p,s} \notin \mathcal{A}_C^{p,s \triangle} \}} \left\{ \theta H_{i+\frac{1}{2}}^{p,s \triangle} + (1-\theta) H_{i+\frac{1}{2}}^{p,s *} \right\}\\
H^-_{i-\frac{1}{2}}& = \sum_{p,s: \{ u_i^{s,p} \in \mathcal{A}_C^{p,s*} \wedge u_i^{p,s} \notin \mathcal{A}_C^{p,s \triangle} \}} \left\{ \theta H_{i-\frac{1}{2}}^{p,s \triangle} + (1-\theta) H_{i-\frac{1}{2}}^{p,s *} \right\}\\
\end{align*}
where:
\begin{equation*}
H^{p,s \triangle}_{i+\frac{1}{2}} =
\begin{cases}
F^{p,s,\triangle}_{i+\frac{1}{2}}  \quad & \mbox{if } F_{i+\frac{1}{2}}^{p,s} \notin \mathcal{A}_F^{p,s*}\\
\\
\dfrac{\Delta t^{p,s}}{\Delta t^{\hat{p},\hat{s}}} F_{i+\frac{1}{2}}^{\hat{p}, \hat{s}\triangle} & \mbox{otherwise}
\end{cases},
\end{equation*}
similarly,
\begin{equation*}
H^{p,s *}_{i+\frac{1}{2}} =
\begin{cases}
F^{p,s,\triangle}_{i+\frac{1}{2}}  \quad & \mbox{if } F_{i+\frac{1}{2}}^{p,s} \notin \mathcal{A}_F^{p,s*}\\
\\
\dfrac{\Delta t^{p,s}}{\Delta t^{\hat{p},\hat{s}}} F_{i+\frac{1}{2}}^{\hat{p}, \hat{s} *} & \mbox{otherwise}
\end{cases},
\end{equation*}
$\hat{s}$ and $\hat{p}$ are superscripts to indicate the previous sub-step of the previous sub-refinement where the flux had been accepted (the last $p,s :  F_{i+\frac{1}{2}}^{p,s} \notin \mathcal{A}_F^{p,s *} \wedge F_{i+\frac{1}{2}}^{p,s} \in \mathcal{A}_F^{p,s \triangle}$).
$\Delta t^{\hat{p},\hat{s}}$, $F_{i+\frac{1}{2}}^{\hat{p}, \hat{s} \triangle}$ and $ F_{i+\frac{1}{2}}^{\hat{p}, \hat{s} *}$ inside the algorithm $S$ are stored in the sets  $\mathcal{T}_F$, and $\mathcal{F}$ so we know their values.
Note that in $p = 0$ we defined $\mathcal{A}_F^* = \emptyset$ so only the first case in the definition of the flux is allowed.

As we said before, the summation at an interface seems to depend on the $i$-th cell that we are considering. It is trivial to show that
$H_{i-\frac{1}{2}}^+ =H_{i-\frac{1}{2}}^-$ if 
for any $ p,s$ such that $$\{ u_i^{s,p} \in \mathcal{A}_C^{p,s*} \wedge u_i^{p,s} \notin \mathcal{A}_C^{p,s \triangle} \} $$ one has that $$\{ u_{i-1}^{s,p} \in \mathcal{A}_C^{p,s*} \wedge u_{i-1}^{p,s} \notin \mathcal{A}_C^{p,s \triangle} \}$$ and vice versa, because in this case both cells have become latent in the same sub-step and the number of evaluated fluxes as their values are the same.
Instead if, for a generic sub-step $\hat{s}$ of a sub-refinement $\hat{p}$ it happens, for example,  that  
$$\{ u_{i-1}^{\hat{s},\hat{p}} \in \mathcal{A}_C^{\hat{p},\hat{s}*} \wedge u_{i-1}^{\hat{p},\hat{s}} \notin \mathcal{A}_C^{\hat{p},\hat{s} \triangle} \}$$
 but $$\{ u_i^{\hat{s},\hat{p}} \in \mathcal{A}_C^{\hat{p},\hat{s}*} \wedge u_i^{\hat{p},\hat{s}} \notin \mathcal{A}_C^{\hat{p},\hat{s} \triangle} \}, $$ this means that the flux $F_{i-\frac{1}{2}}$ has been accepted, because the component $u_{i-1}$ has become latent, but the flux $F_{i+\frac{1}{2}}$ has been rejected in the following sub-step and has to be recomputed, so that a new sub-refinement is required.

The sum at this point can be written for the flux $H^+_{i-\frac{1}{2}}$ as:
\begin{equation*}
\sum_{p,s}^{\hat{p},\hat{s}-1} \left( \theta H_{i-\frac{1}{2}} ^{p,s \triangle} + (1-\theta) H_{i-\frac{1}{2}}^{p,s*} \right) + \theta F_{i-\frac{1}{2}} ^{\hat{p},\hat{s} \triangle } + (1-\theta) F_{i-\frac{1}{2}}^{\hat{p},\hat{s} *}
\end{equation*} 

and for $H^-_{\i-\frac{1}{2}}$ as:
\begin{align*}
\sum_{p,s}^{\hat{p},\hat{s}-1} \left( \theta H_{i-\frac{1}{2}} ^{p,s \triangle} + (1-\theta) H_{i-\frac{1}{2}}^{p,s*} \right) & +  \sum_{n=1}^N \theta \dfrac{\Delta t^n}{\Delta t^{\hat{p}, \hat{s}}} F_{i-\frac{1}{2}} ^{\hat{p}, \hat{s} \triangle} 
+\sum_{n=1}^N  (1-\theta)\dfrac{\Delta t^n}{\Delta t^{\hat{p}, \hat{s}}} F_{i-\frac{1}{2}} ^{\hat{p}, \hat{s} *},
\end{align*}
the $n$ steps are all the later sub-steps of the later sub - refinements where also the flux $F_{i+\frac{1}{2}}$ has been accepted. Due to the recursive nature of the algorithm, we have that $\sum_{n=1}^{N} \Delta t^n = \Delta t^{\hat{p}, \hat{s}}$ because the algorithm exits from the consecutive sub-refinement when the final times are equal,  so that the two different contribution at the end have the same value.
This argument is easily applicable also in the opposite case, when $u_i$ is a latent component while $u_{i-i}$ is an active component. Since
there are no other possible cases, the correct flux balance is preserved at each interface of the domain for each global time steps. 

\subsection{Consistency}
In \cite{hundsdorfer:2007}, explicit multirate schemes for conservation laws have been analyzed, reaching the  conclusion that a method can either be locally inconsistent and mass conservative, or consistent but  not   mass conservative. Here, we will analyse our multirate method in this respect, in the simple case
of the linear advection equation
\begin{equation}
\dfrac{\partial u}{\partial t} + \dfrac{\partial u}{\partial x}=0,
\end{equation}
discretized in space by the finite volume method with a two-point upwind flux.
We assume that at cell $i$ we need to refine the flux $F_{1+\frac{1}{2}}$, while we accept $F_{i-\frac{1}{2}}$. We also assume that 
we perform just one level of refinement by halving the time step. 
If we integrate in time by the forward Euler method, we obtain
\begin{equation}\label{ee_cons}
\begin{split}
u_i^{n+1} &= u_i^{n+\frac{1}{2}} - \frac{1}{\Delta x} \left( F_{i+\frac{1}{2}}^{n+\frac{1}{2}}- \frac{1}{2} {F}_{i-\frac{1}{2}}^n \right)
     	  = u_i^{n+\frac{1}{2}} - \frac{\Delta t}{2 \Delta x} u_i^{n+\frac{1}{2}} + \frac{\Delta t}{2 \Delta x} u_{i-1}^n.
        \end{split}
      \end{equation}

By standard Taylor expansion of the exact solution $u(x,t)$, we have

\begin{align*}
u^{n+1}_i &= u_i^n +\Delta t \frac{\partial u}{\partial t} + \frac{\Delta t^2}{2} \frac{\partial^2 u}{\partial t^2}+  \mbox{h.o.t.}\\ 
u^{n+\frac{1}{2}}_i &= u_i^n +\frac{\Delta t}{2} \frac{\partial u}{\partial t} + \frac{\Delta t^2}{8} \frac{\partial^2 u}{\partial t^2}+ \mbox{h.o.t.}\\
u^{n}_{i-1} &= u_i^n -\Delta x \frac{\partial u}{\partial x}+ \frac{\Delta x^2}{2} \frac{\partial^2 u}{\partial x^2}+ \mbox{h.o.t.}\\
\end{align*}
Replacing into (\ref{ee_cons}), the leading terms of the truncation error $T^{FE}_i$ are 
\begin{equation*}
T^{FE}_i = -\frac{3}{4} \Delta t\frac{\partial^2u}{\partial t^2} -\frac{\Delta t }{2\Delta x} \frac{\partial u}{\partial t} -\frac{\Delta t^2}{8 \Delta x} \frac{\partial ^2 u}{\partial t^2} -\Delta x \frac{\partial^2 u}{\partial x^2}. 
\end{equation*}
As already shown in~\cite{osher:1983}, the truncation error contains the term $\frac{\Delta t }{2 \Delta x} \frac{\partial u}{\partial t}$ which scales as $\frac{\Delta t}{\Delta x}$, and is in general indeterminate for $\Delta t \to 0$ and $\Delta x \to 0. $ Since when studying hyperbolic problems   time and space steps are always reduced  maintaining a constant Courant number, this introduces a consistency error of order $O(1)$.

Instead, if we consider the Backward Euler scheme at time $t^{n+1}$ we get

\begin{equation}\label{ei_cons}
\begin{split}
u_i^{n+1} & = u_i^{n+\frac{1}{2}}-\frac{1}{\Delta x} \left( F_{i+\frac{1}{2}}^{n+1}- \frac{1}{2} {F}_{i-\frac{1}{2}}^{n+1} \right) =u_i^{n+\frac{1}{2}} - \frac{\Delta t}{2 \Delta x} u_i^{n+1} + \frac{\Delta t}{2 \Delta x} u_{i-1}^{n+1}.
\end{split}
\end{equation}
Again, by standard Taylor expansion

\begin{align*}
u^{n+1}_{i-1} &= u_i^n -\Delta x \frac{\partial u}{\partial x} + \Delta t \frac{\partial u}{\partial t} + \frac{\Delta x^2}{2} \frac{\partial^2 u}{\partial x^2} +\frac{\Delta t^2}{2} \frac{\partial^2 u}{\partial t^2}
- 2 \Delta t \Delta x \frac{\partial^2 u}{\partial t \partial x} +  \mbox{h.o.t},
\end{align*}
which plugged into (\ref{ei_cons}), allows to obtain a consistency error $T^{IE}_i$ whose leading terms are

\begin{equation*}
T^{IE}_i= -\frac{3}{4} \Delta t \frac{\partial^2 u}{\partial t^2} + \frac{1}{2} \Delta t \frac{\partial^2 u}{\partial x^2}   -4 \Delta t \frac{\partial^2 u}{\partial t \partial x}.
\end{equation*}

Therefore, consistency is maintained at the expected order.

To generalize these results we consider both explicit and implicit Euler methods where a generic sub-step $(1-\delta) \Delta t$ has been used to go from time $t^{n+\delta}$ to time $t^{n+1}$. 

We use the Taylor expansion centered in a generic point $t^{n+ \tau}$ to obtain
\begin{equation*}
\begin{split}
u_i^{n+1} =  u_i^{n+ \tau} & + (1- \tau) \Delta t \frac{\partial u_i^{n+ \tau}}{\partial t} + \frac{(1 - \tau)^2\Delta t^2 }{2} \frac{\partial^2 u_i^{n+ \tau}}{\partial t^2} +  \mbox{h.o.t.},\\
 u_{i-1}^n =   u_i^{n + \tau} & - \tau \Delta t  \frac{\partial u_i^{n+ \tau}}{\partial t}  + \frac{\tau^2\Delta t^2 }{2} \frac{\partial^2 u_i^{n+ \tau}}{\partial t^2} - \Delta x \frac{\partial u_i^{n+ \tau}}{\partial x}  \\& + \frac{\Delta x^2}{2} \frac{\partial^2 u_i^{n+ \tau}}{\partial x^2} + \Delta x  \tau \Delta t  \frac{\partial^2 u_i^{n+ \tau}}{\partial t \partial x} + \mbox{h.o.t.},
\\ u_i^{n+\delta} = u_i^{n+\tau} & + (\delta-\tau) \Delta t \frac{\partial u_i^{n+\tau}}{\partial t} + \frac{(\delta-\tau)^2}{2} \Delta t^2 \frac{\partial^2 u_i^{n+\tau}}{\partial t^2} + \mbox{h.o.t}.
\end{split}
\end{equation*}    
For the forward Euler we obtain
\begin{equation*}
\begin{split}
T^{FE}_i =&  - \frac{\Delta t}{\Delta x} (\delta-2\tau)\frac{\partial u}{\partial t} - \frac{\Delta t}{2}\left[ \frac{\Delta t}{\Delta x} (\delta^2 - 2\delta\tau)+ (1+\delta-2\tau)\right]\frac{\partial^2 u}{\partial t^2} + \frac{\Delta x}{2} \frac{\partial^2 u}{\partial x^2} +  \tau \Delta t \frac{\partial ^2 u}{ \partial t \partial x}.
\end{split}
\end{equation*}
In this case, the additional term scales as $ \frac{\Delta t}{\Delta x} (\delta-2\tau)$.

The implicit Euler method is instead consistent for any value of  $\tau$ since in this case we have
\begin{equation*}
\begin{split}
u_{i-1}^{n+1} =   u_i^{n + \tau} & +(1- \tau) \Delta t  \frac{\partial u_i^{n+ \tau}}{\partial t}  + \frac{(1-\tau)^2\Delta t^2 }{2} \frac{\partial^2 u_i^{n+ \tau}}{\partial t^2} - \Delta x \frac{\partial u_i^{n+ \tau}}{\partial x} + \frac{\Delta x^2}{2} \frac{\partial^2 u_i^{n+ \tau}}{\partial x^2} \\ &- \Delta x  (1-\tau) \Delta t  \frac{\partial^2 u_i^{n+ \tau}}{\partial t \partial x} + \mbox{h.o.t.},
\end{split}
\end{equation*}
and thus,
\begin{equation}
\begin{split}
T^{IE}_i =& -\frac{\Delta t}{2 } (1-\delta-2\tau)\frac{\partial^2 u}{\partial t^2}+ \frac{\Delta x}{2} \frac{\partial^2 u}{\partial x^2}  + \Delta t (1-\tau) \frac{\partial ^2 u}{ \partial t \partial x} 
\end{split}
\label{taylor_gen_impl}
\end{equation}

From these considerations we can deduce that, if we use the $\theta$ method, we would have an inconsistent scheme whenever $\theta \neq 1$. The inconsistency term is also present for the TR-BDF2 scheme that we introduce in the next Section.  Therefore, any   conservative multirate scheme not based on the backward Euler
method would introduces a consistency error analogous to that discussed in \cite{hundsdorfer:2007}
for explicit schemes. It can be argued, however, that
 this fact does not reduce the effectiveness of such methods for practical applications. Indeed, the goal of a multirate approach is to reduce the computational cost by using a relatively large $\Delta t$ and
 refining  it  only  in the region where is necessary to keep the discretization error small. The error is controlled by setting the appropriate tolerance in the algorithm $S$ which accept/reject the fluxes for a given space discretization.
A problem may however arise if the multirate scheme is combined with dynamic adaption in space. For this situation the effect of the consistency error in conservative multirate schemes has to be investigated further, but this is beyond the scope of the present work. 

\subsection{Time discretization with TR-BDF2}
While a time discretization based on the $\theta-$method has been employed to introduce the proposed
conservative multirate method and for the consistency analysis, for the numerical experiments and the practical application of the present
approach we have exploited, as in~\cite{bonaventura:2018}, the 
 TR-BDF2 method, because of its interesting properties.  This method is a composite one step, two stages method, consisting of one stage of the trapezoidal scheme followed by one stage of the BDF2 method.
It can be written for the discretization of an ODE system
$y^{\prime}=f(t,y) $ as
\begin{align*}
& u^{n + \gamma} = u^n + \frac{\Delta t_n \gamma}{2} \left( f(t_n,u^n) + f(t_{n+ \gamma},u^{n+\gamma}) \right) \\
& u ^ {n+1} = \frac{1}{\gamma (2- \gamma)} u^{n + \gamma} - \frac{(1-\gamma)^2}{ \gamma (2- \gamma)} u^{n} + \frac{1-\gamma}{2-\gamma} \Delta t_n f(t_{n+1} , u^{n+1})
\end{align*} 
For $\gamma = 2 - \sqrt{2}, $ the method is L-stable and  also employs the same Jacobian matrix for the two stages.   
In~\cite{hosea:1996} it has been interpreted as a Diagonally Implicit Runge Kutta (DIRK) method with two internal stages, proving the following properties:
\begin{itemize}
\item the method is strongly S-Stable;
\item it is endowed with a Cubic Hermite interpolation algorithm that yields globally $\mathcal{C}^1$ continuous trajectories.
\end{itemize} 
Due to its favorable properties, it has been recently applied for efficient
discretization of high order finite element methods for numerical weather forecasting in~\cite{tumolo:2015},
while its monotonicity properties have been studied in~\cite{bonaventura:2017}.

\subsection{Flux-partitioning and error estimator}
To select the components that have to be recomputed with a smaller time step,  we need to introduce a local error estimator for the fluxes.
A simple approach is to compare the fluxes computed with the $\theta$-method or the TR-BDF2 method, with the fluxes at the same interface cell computed with a more accurate method. The absolute value of the difference between the two fluxes can be used as a measure of the error.
For $\gamma=2-\sqrt{2}$ the TR-BDF2 scheme has a third order method embedded, this fact can be exploited to derive the error estimator, yet
as remarked in~\cite{hosea:1996}, the third order method embedded in TR-BDF2 is not
A-stable. In that work a heuristic approach that entails the solution of an additional linear system per time step
has been proposed to stabilize the error estimator. For a large ODE systems coming from the spatial discretization of PDEs,
solving at each time step this extra linear system could turn out to be very expensive. 

Therefore, we propose another types of error estimator, which are less
expensive.  At each time step, for a two stage method as the TR-BDF2 method, we know the
active components values at times $t_n$ and $t_{n+\gamma}$ , so we can use an extrapolation technique to obtain a prediction of the value at time $t_{n+1}$.
If we call the extrapolated solution at time $t_{n+1}$ as $\bar{u}_{ext}^{n+1}$, the extrapolated fluxes at the interface are $\bar{F}_{ext_{i+\frac{1}{2}}}^{n+1}$ and we obtain the error estimator as:
\begin{equation*}
\mathcal{R}_F = \{ F_{i+\frac{1}{2}} : |F_{i+\frac{1}{2}} ^{n+1}-\bar{F}_{ext_{i+\frac{1}{2}}}^{n+1}| > \tau_r|F_{i+\frac{1}{2}} ^{n+1}| + \tau_a \}
\end{equation*}

The simplest extrapolation technique is the linear extrapolation, given by
\begin{equation*}
\bar{u}_{lin}^{n+1} = u^n + \dfrac{t_{n+1}-t_{n}}{t_{n+\gamma}-t_{n}}\left( u^{n+\gamma} - u^{n} \right),
\end{equation*}
by which we obtain the extrapolated values of  $F_{lin_{i+\frac{1}{2}}}^{n+1}$ at the required interface, whose difference with the computed value provides the error estimator.

A more precise estimator can be obtained by applying a cubic Hermite extrapolation at time $t^n$ and $t^{n+\gamma}$ considering the fact that the TR-BDF2 method provides a formula to compute the coefficient for the cubic Hermite extrapolation easily.

The extrapolation can be evaluated as:
\begin{eqnarray}
\bar{{u}}_{cub}(t) &=& ( \alpha_3 - 2 \alpha_2)\beta(t)^3 
 + (3\alpha_2 - \alpha_3)\beta(t)^2   
  + \alpha_1\beta(t) +\alpha_0,  \nonumber
 \label{hermite}
\end{eqnarray}
	
$\alpha$ coefficients are:
	\begin{equation*}
\begin{aligned}
	\alpha_0&=u^n, \hspace{1em} \alpha_1=\gamma \Delta t_n f(t_n,u^n), \hspace{1em} \alpha_2=u^{n+\gamma}-u^n - \alpha_1, \nonumber \\
	\alpha_3&= \gamma \Delta t_n ( f(t_{n+\gamma}, u^{n+\gamma}) -  f(t_{n}, u^{n}) ), \nonumber
	 \end{aligned}
\end{equation*}
instead $\beta$ is:
$$
 \beta(t) = \frac{t-t_n}{\gamma \Delta t_n}.
 $$
 At time $t^{n+1}$ the extrapolated solution would be:
 $$ 
\bar{{u}}_{cub}^{t+1} = ( \alpha_3 - 2 \alpha_2)\left( \frac{1}{\gamma}\right)^3 
 + (3\alpha_2 - \alpha_3)\left( \frac{1}{\gamma}\right)^2   
  + \alpha_1\left( \frac{1}{\gamma}\right) +\alpha_0,   
$$ 
 
In our test cases we use the error estimator based on the Cubic Hermite extrapolation.

\subsection{Systems of PDEs}
\label{system_ref}
The multirate method is easily extended to a system of non-linear conservation laws. The only non trivial part is how to define the set of active fluxes.   

A system of $d$ non-linear conservation laws can be written as:
\begin{equation}
\label{systemPDE}
\frac{\partial {\bf u}}{ \partial t} + \frac{\partial {\bf(f(u))}}{\partial x} = 0 \qquad x \in \mathbb{R}\quad t>0 
\end{equation}
where $\bf u$ and $\bf f$ are $d$-vectors on the problem domain, ${\bf u} = [u_1, u_2, \cdots, u_d]^T$ and 
$${\bf{F(u)}} = [F_1(u_1,\cdots, u_d), F_2(u_1,\cdots, u_d), \cdots, F_d(u_1,\cdots, u_d)]^T$$
is a vector of fluxes.

If we use a two-point flux approximation, when (\ref{systemPDE}) is semi-discretized in space, the flux at each interface depends on the values at the right and at the left cell of all variables $u_1, \cdots, u_d$.
To preserve the mass of the whole system, if the $j$-th flux for the $i$-th variable has been rejected by our error estimator, all fluxes at the same space position have to be considered as rejected.   

\begin{figure}
\includegraphics[scale=0.9]{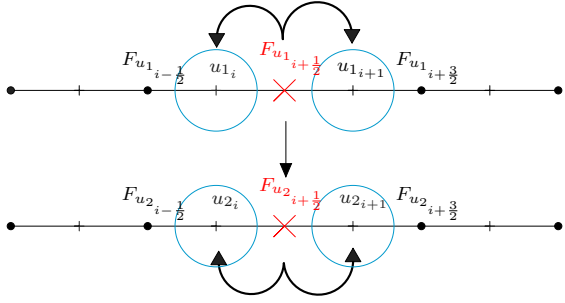}
\caption{Example of rejected fluxes in a system of non-linear conservation laws}
\label{ex:fluxsystem}
\end{figure}

In Fig. \ref{ex:fluxsystem}, we show a simple example with $d=2$. If the flux for the variable $u_1$ has been rejected in position $x_{i+\frac{1}{2}}$, the components $u_{1_i}$ and $u_{1_{i+1}}$ will be included in the set of active components but, to be conservative, also the flux for the variable $u_2$ will be rejected and so also the components  $u_{2_i}$ and $u_{2_{i+1}}$ will be recomputed with a smaller time step.
 
\subsection{Boundary conditions}
To illustrate our scheme we have assumed that the differential problem is set on the whole real line. However, in the numerical tests of the next Section (as well as in all practical situations) we have to deal with bounded domain, and proper boundary conditions must be imposed. Since we are adopting a finite volume scheme, the boundary conditions have been applied by computing the fluxes at the fictitious boundary interface by  the well known ``ghost node'' technique.
With this method the correct type of information (i.e. that corresponding to the characteristics entering the domain) is automatically selected by the numerical scheme.

\section{Numerical experiments}
In this section, we present different numerical experiments to test the efficiency and the accuracy of the conservative multirate method. First we show the multirate method applied to the Burgers' equation, then a more complex scalar test case, the Buckley-Leverett equation and, at the end, we illustrate the multirate method applied to a system of nonlinear conservation laws, the Shallow Water equations.

\label{numerical}
\subsection{Burgers equation}
Here, we apply the multirate method to Burgers equation with Dirichlet boundary conditions, thus
repeating the tests presented in \cite{bonaventura:2018}, but with the conservative variant of our algorithm.
The Burgers equation is a nonlinear conservation law given by

\begin{equation*}
\begin{cases}
\dfrac{\partial u}{\partial t}+ \dfrac{\partial}{\partial x}\left(\dfrac{1}{2}u^2\right)=0 & (x,t) \in (-1,3) \times (0,1),\\
u(x,0)=u_0(x) & x\in(-1,3),\\
u(-1,t) = u_l(t) \quad u(3,t)=u_r(t) &t\in(0,1),
\end{cases}
\end{equation*}
where $
u_0(x) = \begin{cases}
u_l(t) \qquad &x<0,\\
u_r(t) \qquad &x>0.
\end{cases}
$

The form of the solution depends on the relation between $u_l$ and $u_r$. 

\subsection*{First case: $u_l > u_r$} 
\begin{figure*}
\begin{center}
\includegraphics[width=0.3\textwidth]{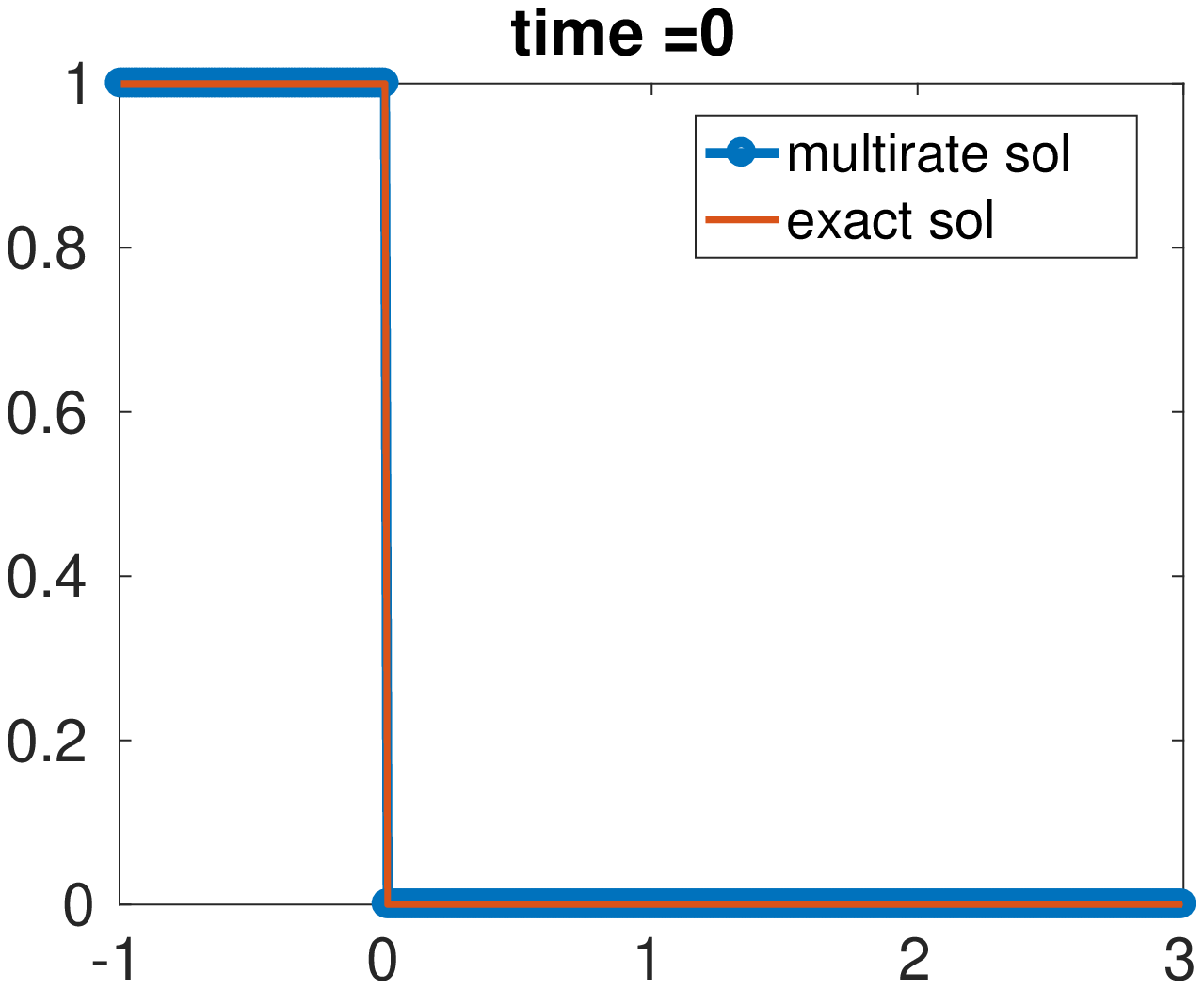}
\includegraphics[width=0.3\textwidth]{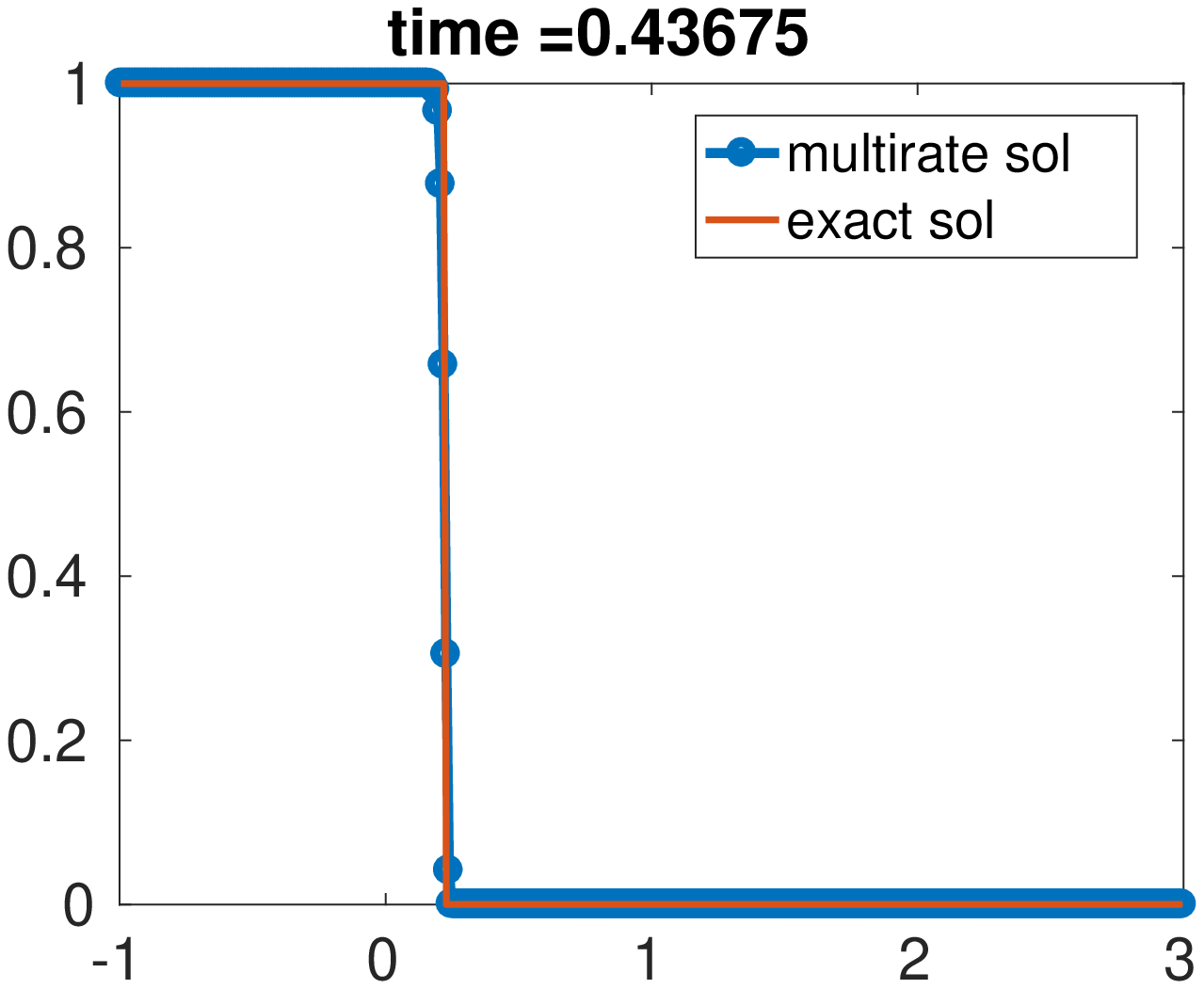}
\includegraphics[width=0.3\textwidth]{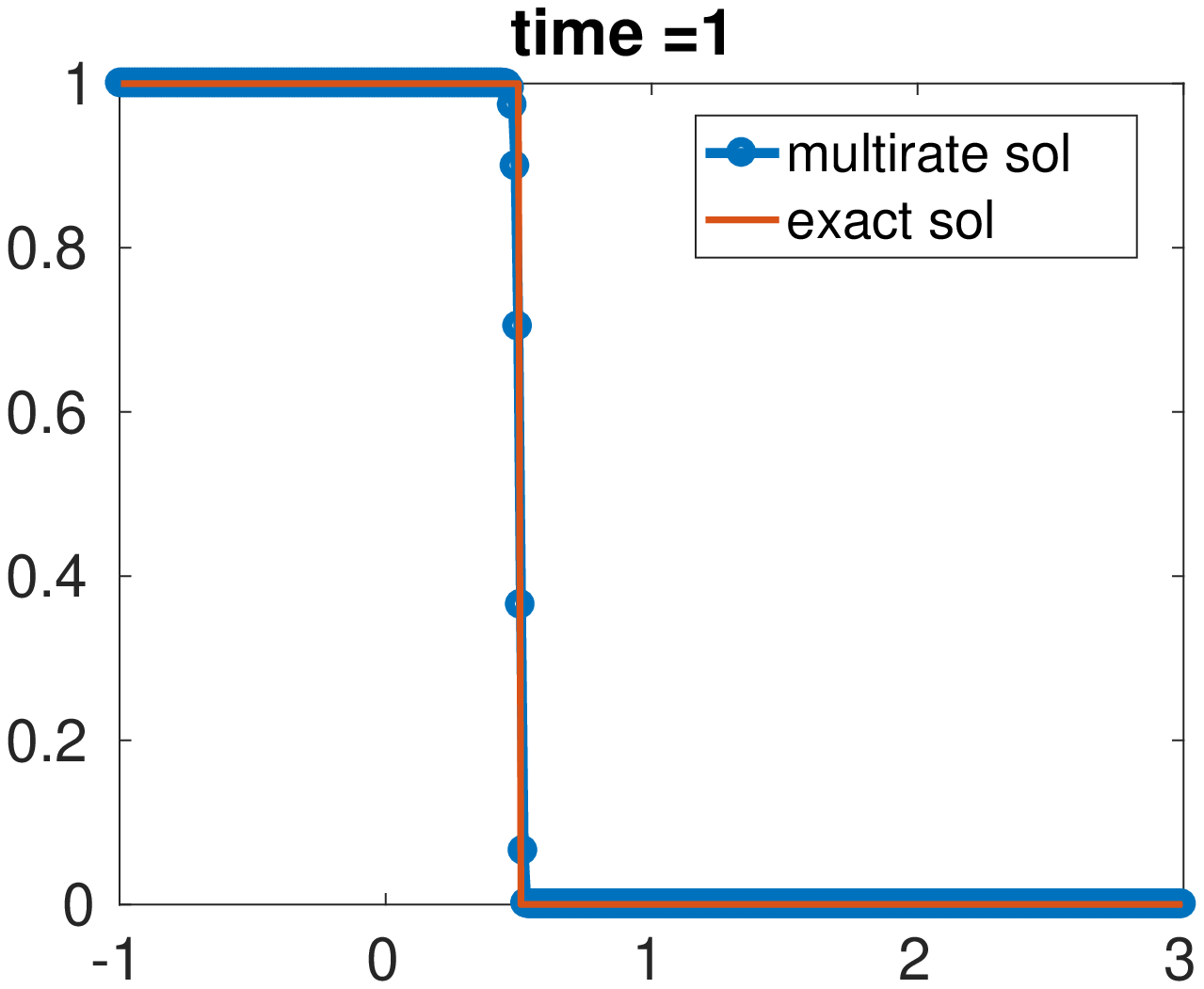}
\end{center}
\caption{Multirate TR-BDF2 integration and exact solution for the shock wave at different times $ t=0$s, $ 0.45$s and $1$s.}
\label{UBurgers1}
\end{figure*}

\begin{figure}
\centering   
\includegraphics[width=0.6\textwidth]{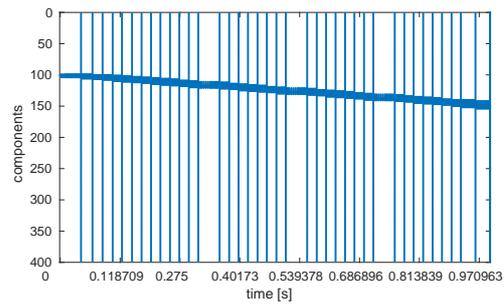}
\caption{The components being computed at each time step by the TR-BDF2 method for the burgers equation that generates a shock wave.}
\label{Burgers1mesh}
\end{figure}

\begin{figure}
\centering 
\includegraphics[width=0.6\textwidth]{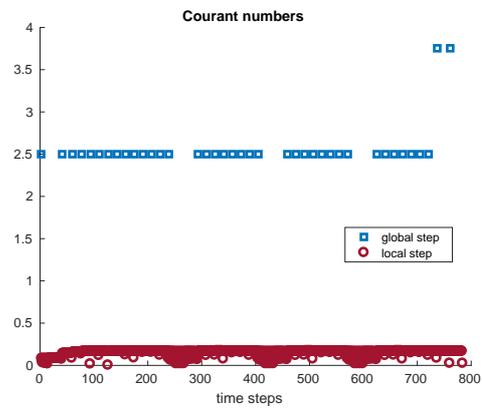}
\caption{Courant number for each time step for the shock wave.}
\label{Burgers1Courant}
\end{figure}

In this case we consider $u_l =1$ and $u_r = 0$
 with a number of cells equal to $400$, the absolute and relative error tolerances are $ 10^{-4}, $ $ 10^{-6},$
 respectively, while the tolerance for the Newton solver is $10^{-14}$ on the difference between two consecutive iterations.
The TR-BDF2 method has been used as solver to integrate in time, the size of the global time step is equal to $0.1s$. To obtain an entropic solution we used the local Lax Friedrichs flux~\cite{toro:2013} (also know as Rusanov flux) as numerical flux for the two point Finite Volume method:
\begin{equation}
F_{i+\frac{1}{2}}=F_{i+\frac{1}{2}}(u_i,u_{i+1}) = \frac{1}{2} \left[ (f(u_{i+1})+ f(u_i)) - \alpha (u_{i+1} - u_{i}) \right],
\label{rusanov}
\end{equation}
where $\alpha = \max_\omega  |f'(\omega)|$ 
and the maximum is taken in the range $\omega \in \left[ u_i, u_{i+i} \right]$.
As we can see in Fig. \ref{UBurgers1}, the solution computed with the multirate method is in excellent agreement with as the exact solution.  
 In Fig. \ref{Burgers1mesh} we represent the set of active components at each time. We can observe that the multirate method captures the shock and refines only the region of the domain where the solution is changing rapidly. We also plot the Courant numbers for each time step, Fig. \ref{Burgers1Courant}. The self adjusting strategy selects small Courant numbers inside the time slab, while the global step correponds to a Courant number equal to $2.5$. Note that we prescribed a global step size equal to $0.1,$ that gives a Courant number of $10$, but all components have been rejected for the given value of the error tolerance,  so that the global time step size is in fact smaller and equal to $0.025$s except for the last two time slabs.

\subsection*{Second case: $u_r > u_l$}
To obtain a rarefaction wave, we set the value at the left $u_l=0$ and the value at the right $u_r=1$. The boundary conditions are $u(-1,t)=u_{l} \quad \forall t \in (0,1)$ and $u(3,t)= u_r \quad \forall t \in (0,1),$ while the other parameters are the same as in the previous test case.

\begin{figure*}
\begin{center}
\includegraphics[width=0.3\textwidth]{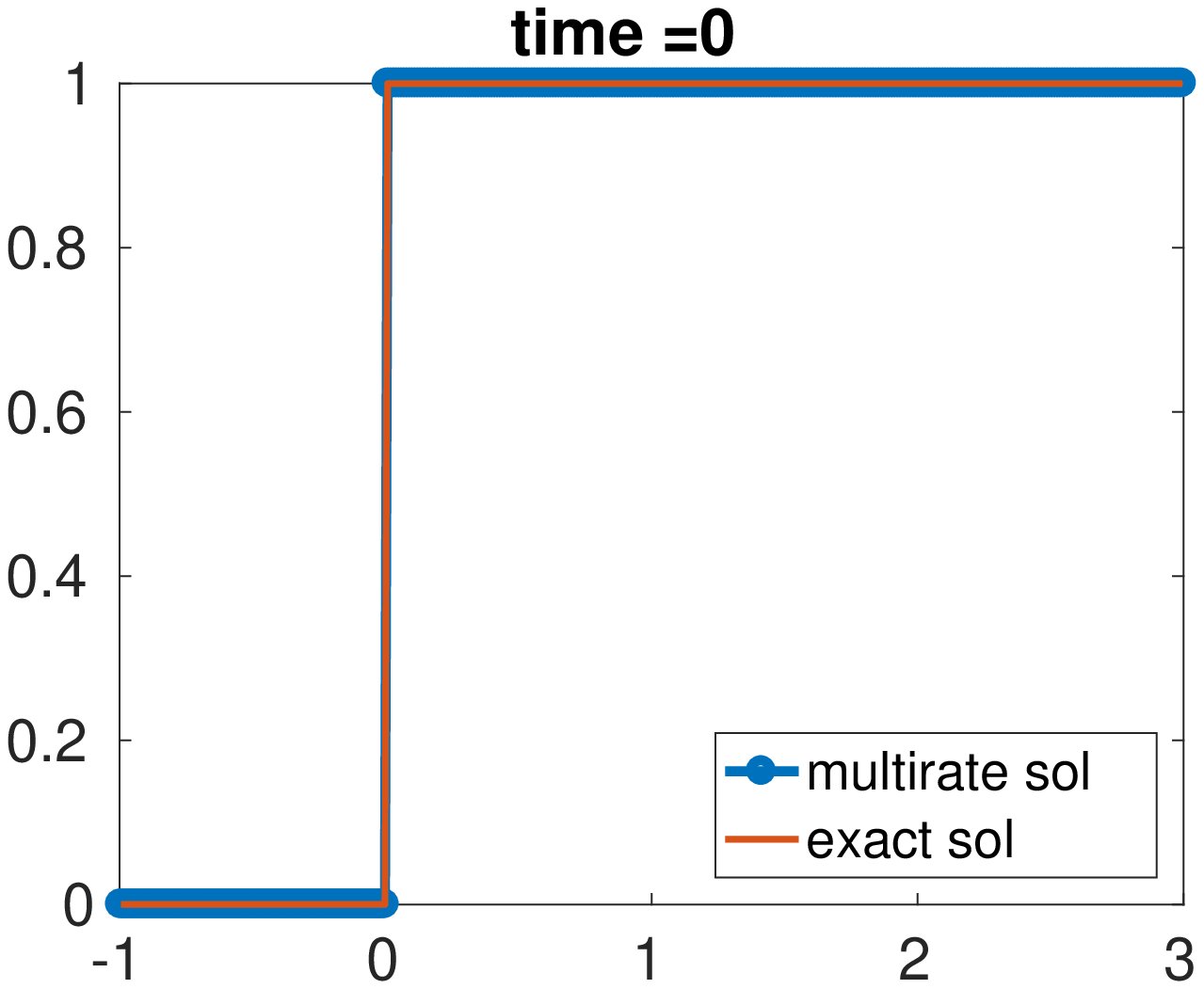}
\includegraphics[width=0.3\textwidth]{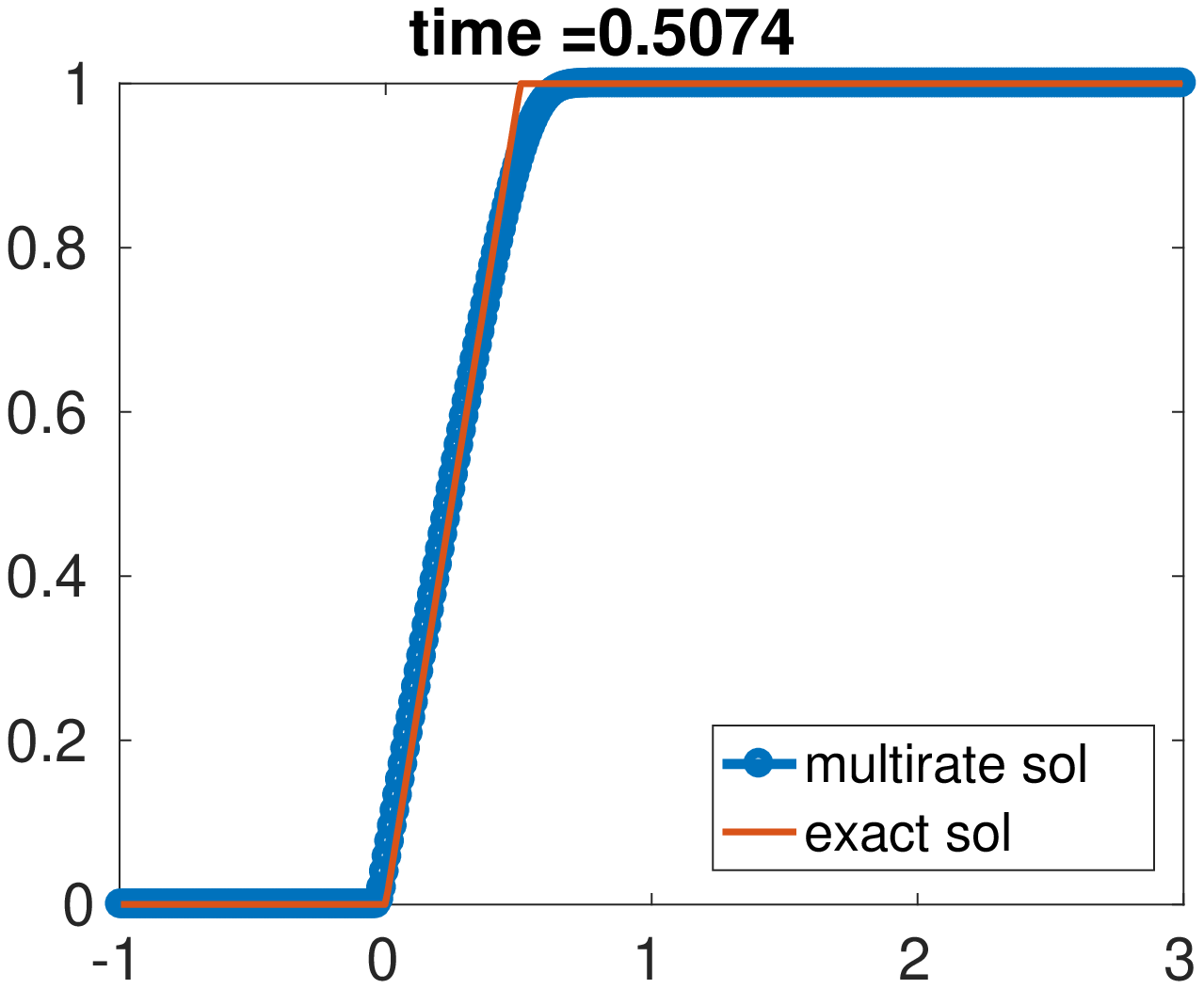}
\includegraphics[width=0.3\textwidth]{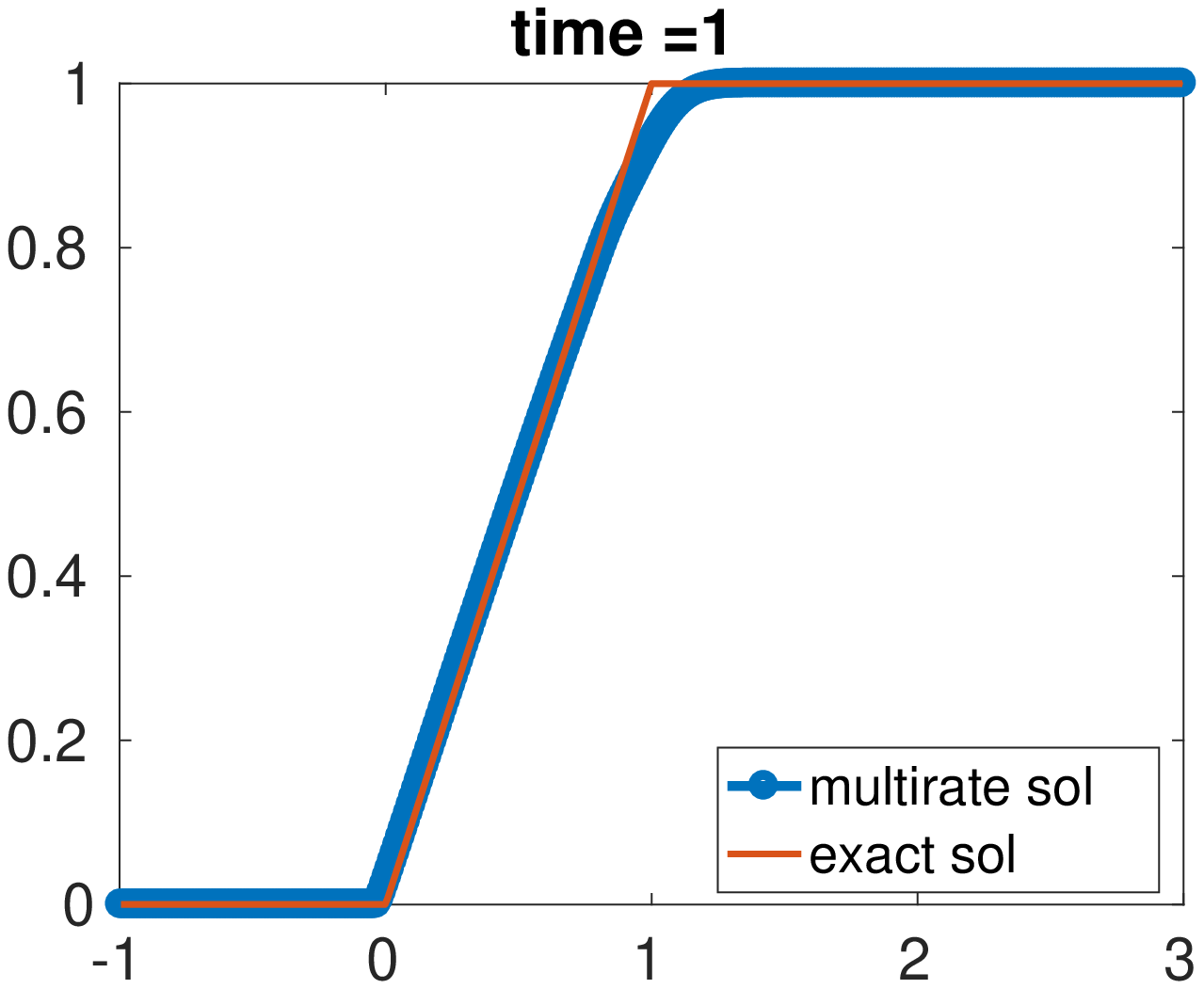}
\end{center}
\caption{Multirate TR-BDF2 integration and the exact solution for the rarefaction wave at different times.}
\label{UBurgers2}
\end{figure*}
\begin{figure}
\centering   
\includegraphics[width=0.6\textwidth]{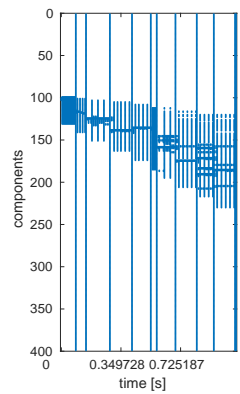}
\caption{The components being computed at each time step with the TR-BDF2 method for the burgers equation that generates a rarefaction wave.}
\label{Burgers2mesh}
\end{figure}
\begin{figure}
\centering
\includegraphics[width=0.6\textwidth]{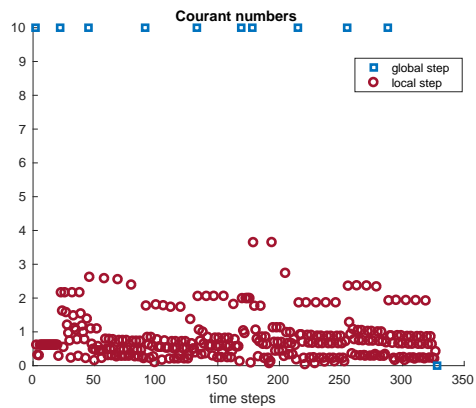}
\caption{Courant number for each time step for the rarefaction wave.}
\label{Burgers2Courant}
\end{figure}
In Fig. \ref{UBurgers2} we can see   the solution obtained with the multirate method. The
numerical diffusion is clearly visible due to the first order monotone flux employed. 
In this case, the Courant number for the global step is equal to $10$, as shown in Fig. \ref{Burgers2Courant}. The Courant numbers for the step inside the time slab are larger than those obtained in the shock wave solution and less time steps are necessary to compute the solution   at the final time.
Fig. \ref{Burgers2mesh} represents the set of active components at each time. As expected, the size of the set increases with time because the rarefaction zone is expanding.

\subsection{Buckley-Leverett equation}
An example of a more complex conservation law is given by the Buckley-Leverett equation:
\begin{equation*}
\begin{cases}
\dfrac{\partial u}{\partial t} + \dfrac{\partial}{\partial x}f(u) = 0 & (x,t) \in (0,2\pi)\times(0,1)\\
f(u) = \dfrac{u^2}{u^2+\frac{1}{3}(1-u)^2}\\
u(x,0) = \sin(x) & x \in (0, 2\pi)\\
u(0,t)= u(2\pi,t) & t \in (0,1)
\end{cases}
\end{equation*}
Also in this case, we used two-point finite volumes with Rusanov flux, with  $N_x = 100$  cells. To integrate up to time $T = 0.5$ the TR-BDF2 method has been used with  a global size step $\Delta t = 0.1$. In this case,   periodic boundary conditions were employed. The absolute and relative error tolerances are $ 10^{-4}, $ $ 10^{-5},$  respectively, while the tolerance for the Newton solver is $10^{-13}$.  To compute the l1-norm of the error we use as a refenrece solution that provided by the Matlab solver \texttt{ode45} with maximum time step allowed equal to $\Delta t = 10^{-5}$s.

\begin{figure*}
\begin{center}
\includegraphics[width=0.3\textwidth]{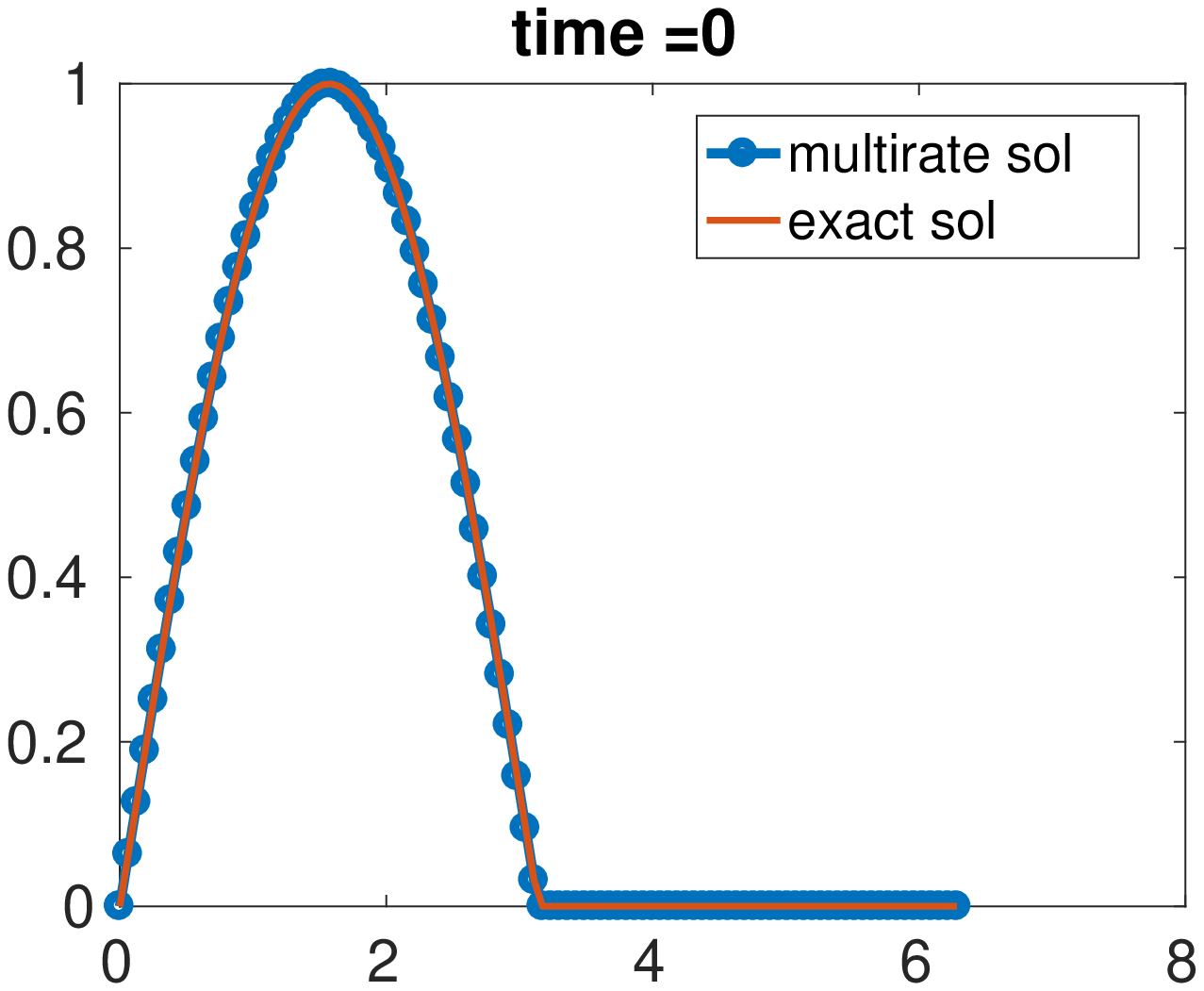}
\includegraphics[width=0.3\textwidth]{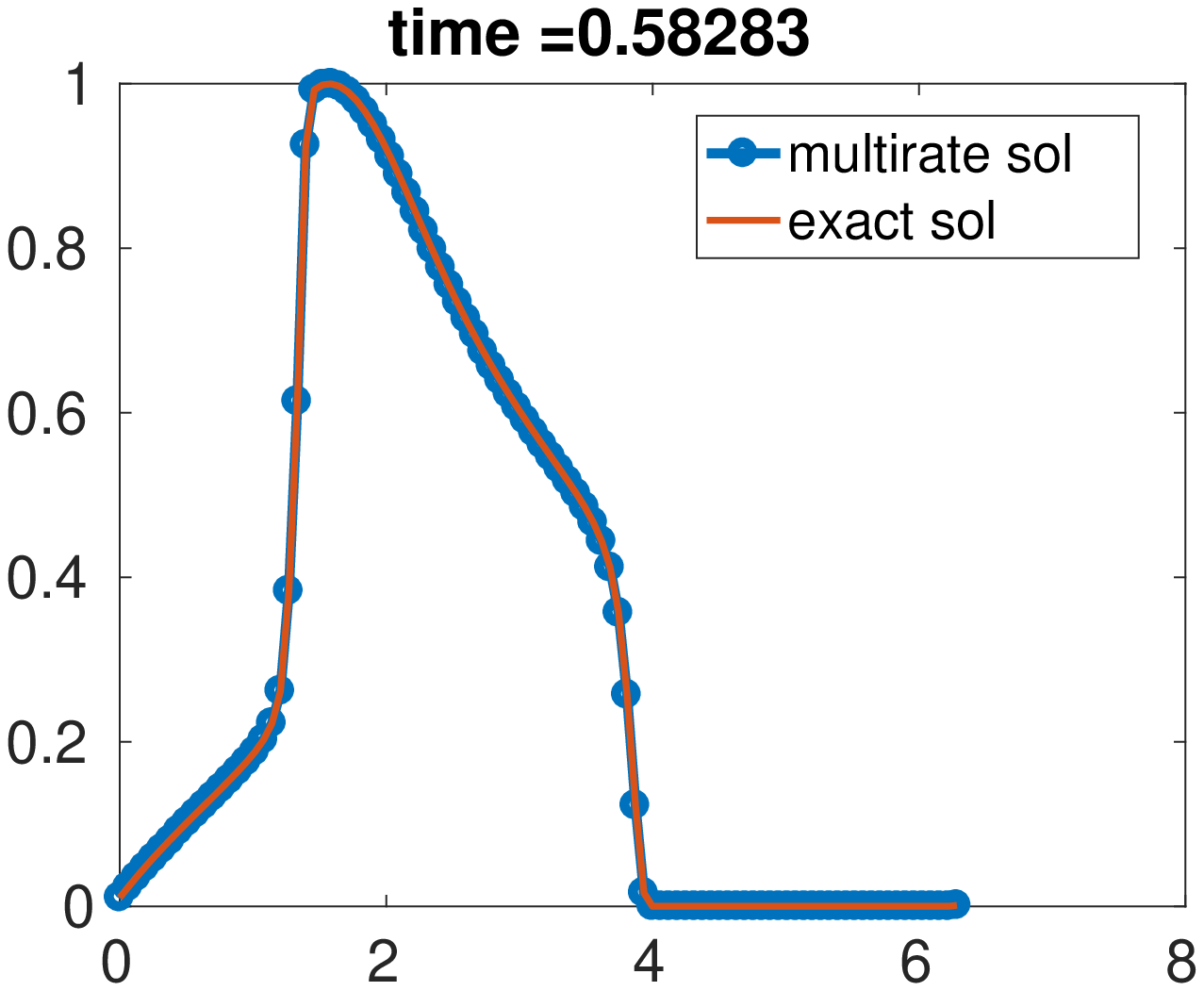}
\includegraphics[width=0.3\textwidth]{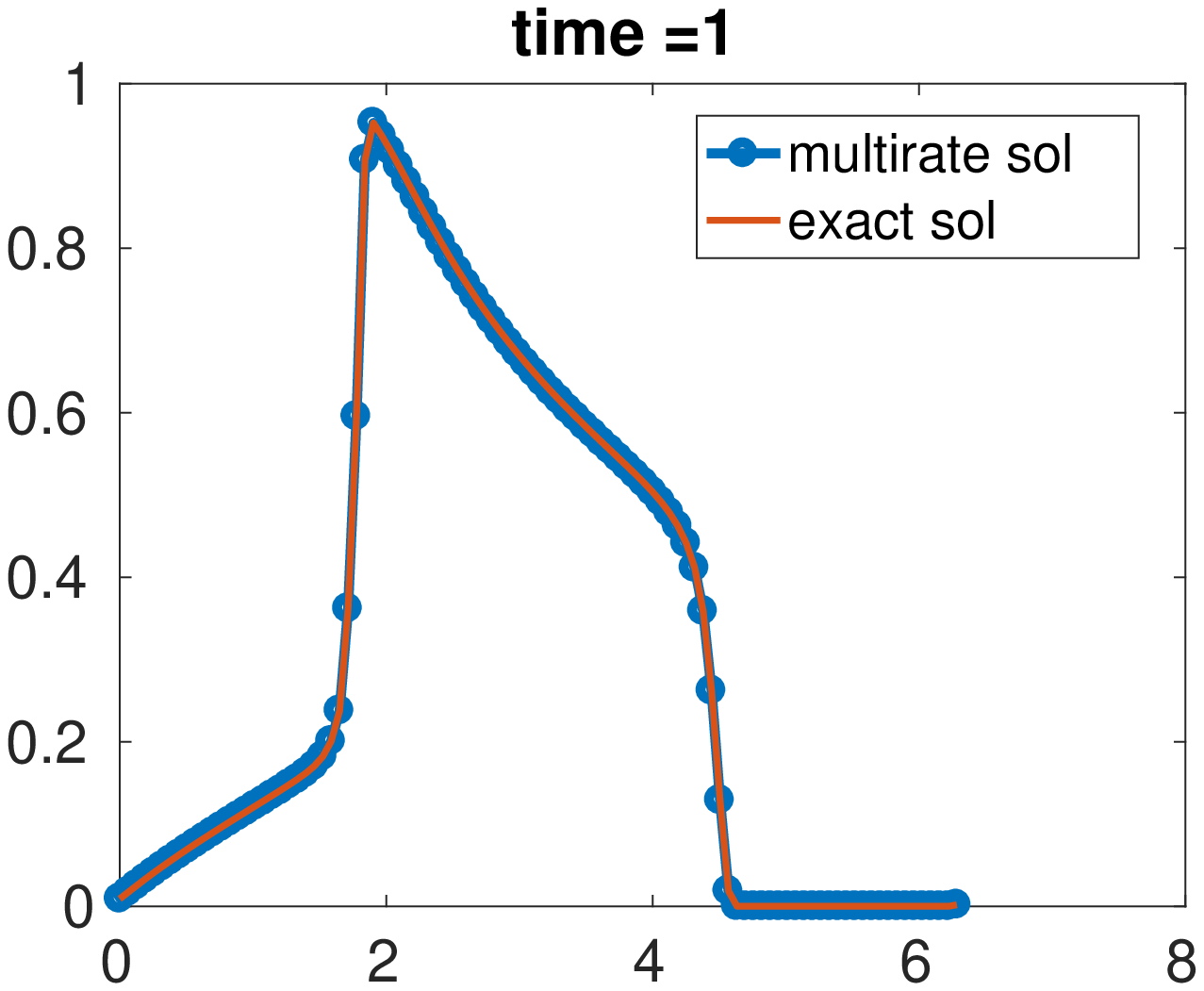}
\end{center}
\caption{Multirate TR-BDF2 solution and the solution computed with the \texttt{ode45} matlab solver.}
\label{bl:sol}
\end{figure*}

\begin{figure}
\centering
\hspace*{-2cm}
\includegraphics[width=0.9\textwidth]{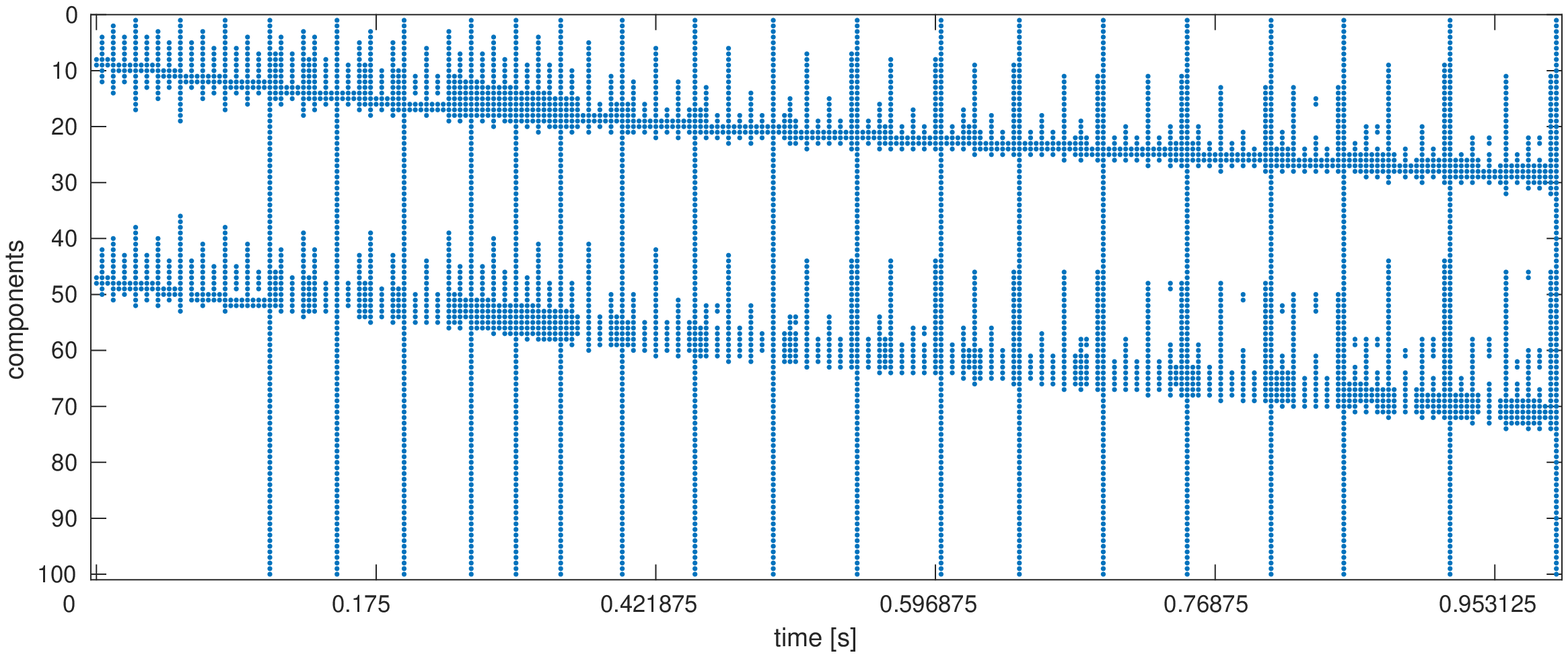}
\caption{The components beig computed at each time step by the TR-BDF2 method for the Buckley-Leverett problem.}
\label{BL:mesh}
\end{figure}

\begin{figure}
\centering
\includegraphics[width=0.6\textwidth]{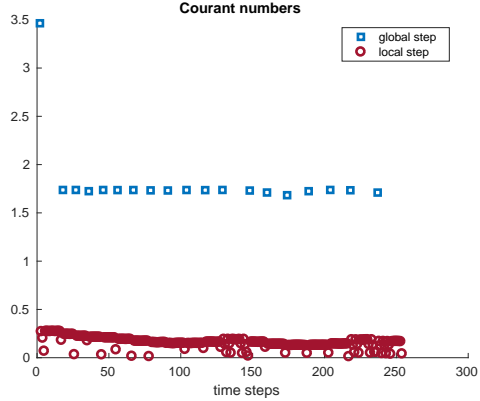}
\caption{Courant number for each time step for the Buckley-Leverett equation.}
\label{BL:Courant}
\end{figure}

This is a more complex test case, because of both a shock and a rarefaction wave appear in the solution, as we can see in Fig.~\ref{bl:sol}. The multirate method refines the solution only where the solution is moving very fast (Fig. \ref{BL:mesh}) using smaller Courant numbers, as illustrated in~Fig.~\ref{BL:Courant}.   

We then compare our mass conservative approach with the original multirate method proposed in~\cite{bonaventura:2018}. As shown in Table \ref{BL:table}, we obtain essentially the same error in the $l_1$-norm for both  methods, but, while with the previous method the system loses $4 \%$ of the mass during the simulation, the new method, preserves the total mass of the system as expected. 

\begin{table}
\caption{The ratio between the initial and final mass, the normalized difference between the initial and final mass
in the Buckley-Leverett equation test case.}
\label{BL:table}
\centering
\begin{tabular}{llll}
\hline\noalign{\smallskip}
 & ratio  & diff. & $l1$-norm \\
\noalign{\smallskip}\hline\noalign{\smallskip}
MC scheme		  & $1$ & $8.36e-15$ &$0.0013$\\
\noalign{\smallskip}\hline\noalign{\smallskip}
N-MC scheme & $0.96$ & $0.0313$& $0.0012$ \\
\noalign{\smallskip}\hline
\end{tabular}
\end{table}

\subsection{Saint-Venant equations: dam break problem}
We have applied our multirate strategy to the Saint-Venant (or shallow water)
equations, which can be written in   conservative form as:
\begin{equation*}
\begin{cases}
\frac{\partial h}{\partial t}+ \frac{\partial q}{\partial x}=0\\
\frac{\partial q}{\partial t}+ \frac{\partial}{\partial x} \left( \frac{q^2}{h} + g \frac{h^2}{2}
\right) =0.
\end{cases}
\end{equation*}
Here, $h$ denotes the fluid depth and $q=hu $ the discharge, where $u$ is the  velocity of the fluid.
These equations are the core of many numerical models for river hydraulics and environmental
flows. A more complete discussion of the Saint-Venant equations  
can be found in~\cite{leveque:1992}. It has to be remarked that even very efficient single rate semi-implicit
methods, see e.g.~\cite{rosatti:2011}, when applied to the Saint-Venant equations in presence of shocks,
must employ relatively small time steps throughout the domain. As we will see, this shortcoming
is overcome by our approach.

The dam break problem is a special case of the Riemann problem, where at the initial time  $h_0(x) = \begin{cases}
h_l \ & \mbox{if } x<x_0 \\
h_r \ & \mbox{if } x>x_0
\end{cases}$
 and $u=q=0$  everywhere in the domain.
For the spatial discretization of the Saint-Venant equations we used again the Rusanov flux. In this case, the
numerical diffusion coefficient $\alpha$ in (\ref{rusanov}) is defined as:
\begin{equation*}
\alpha = max \{ |\lambda_i^1|, |\lambda_i^2|,|\lambda_{i+1}^1|, |\lambda_{i+1}^2| \},
\end{equation*}
$\lambda^1_i$ and $\lambda^2_i$ are eigenvalues of the system for the control volume $i$:
\begin{align*}
&\lambda_i^1 = \frac{h_i}{q_i} - \sqrt{gh_i}\\
&\lambda_i^2 = \frac{h_i}{q_i} + \sqrt{gh_i}.
\end{align*}

We used  300 cells over the domain $[0,3000], $ while the absolute and relative error
 tolerances are $ 10^{-2}, $ $ 10^{-4},$
 respectively, while the tolerance for the Newton solver is $10^{-13}$.
The size of the global steps is equal to $8$s,  and we integrate in the time interval $[0, 100].$
The initial condition for the water height is $h_0(x)=\begin{cases}
1.5 \ & \mbox{if } x<1500 \\
0 \ & \mbox{if } x>1500
\end{cases}$
and for water velocity $ u = \frac{q}{h}=0$.   

\begin{figure}
\centering
\includegraphics[width=0.4\textwidth]{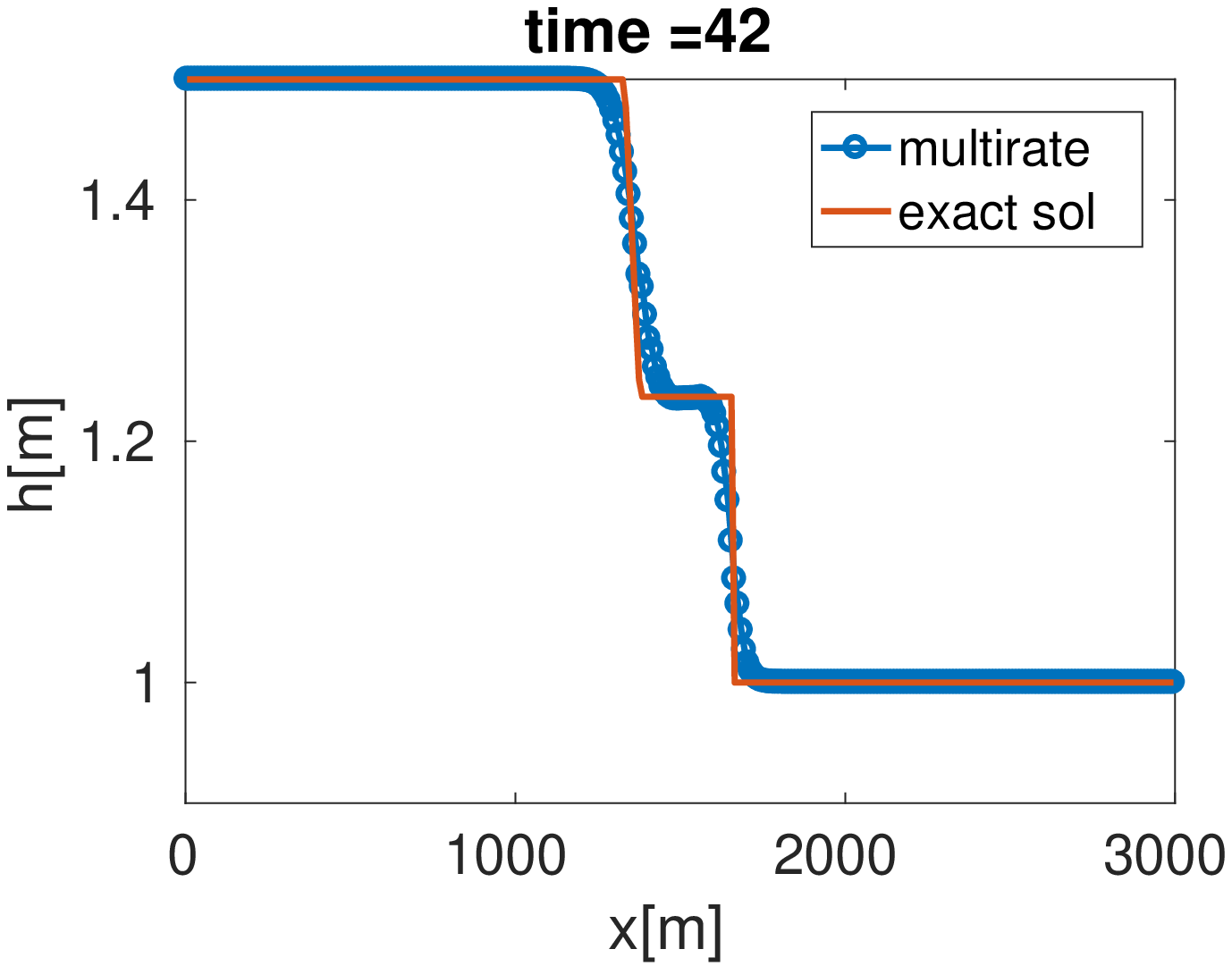}
\includegraphics[width=0.4\textwidth]{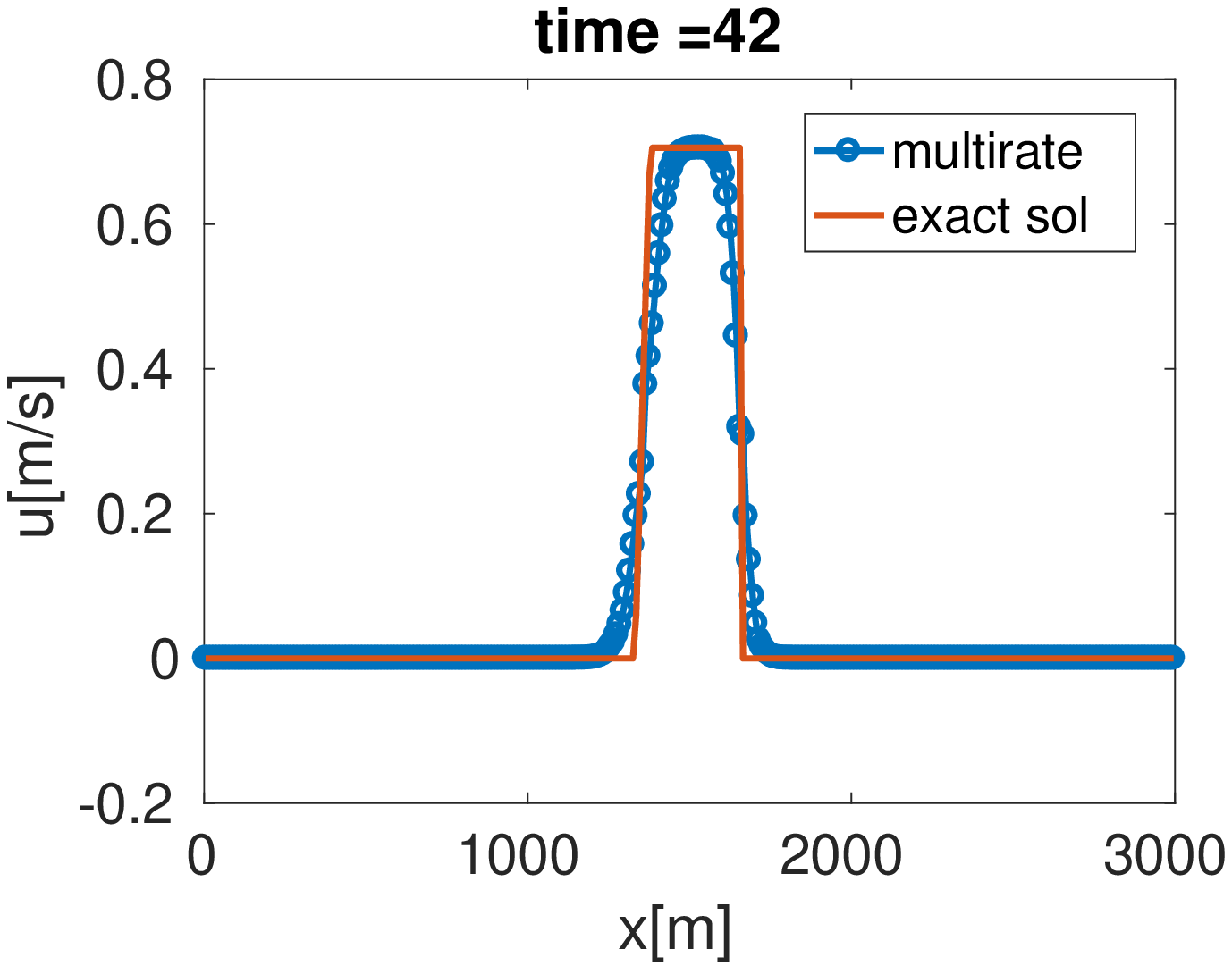}
\includegraphics[width=0.4\textwidth]{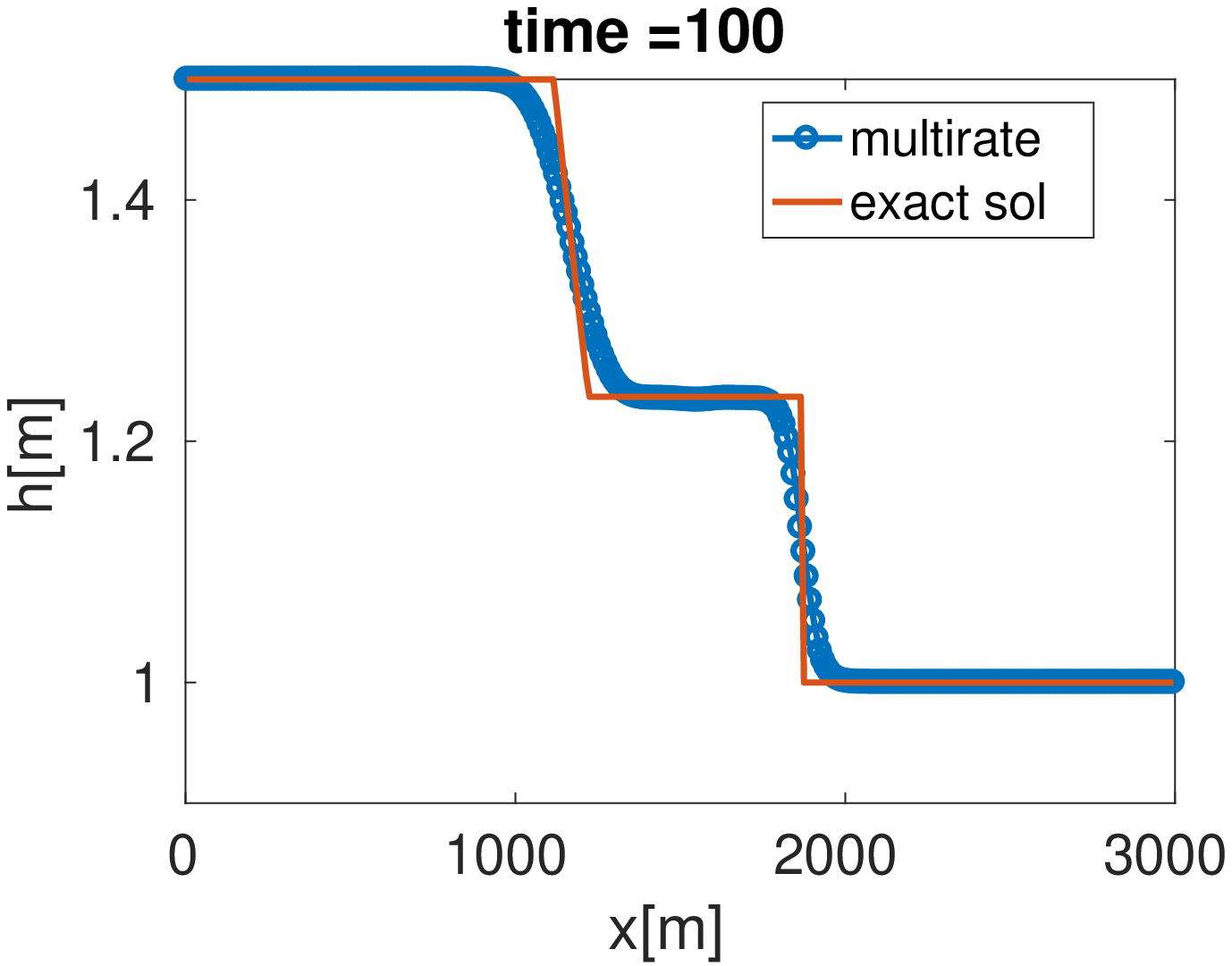}
\includegraphics[width=0.4\textwidth]{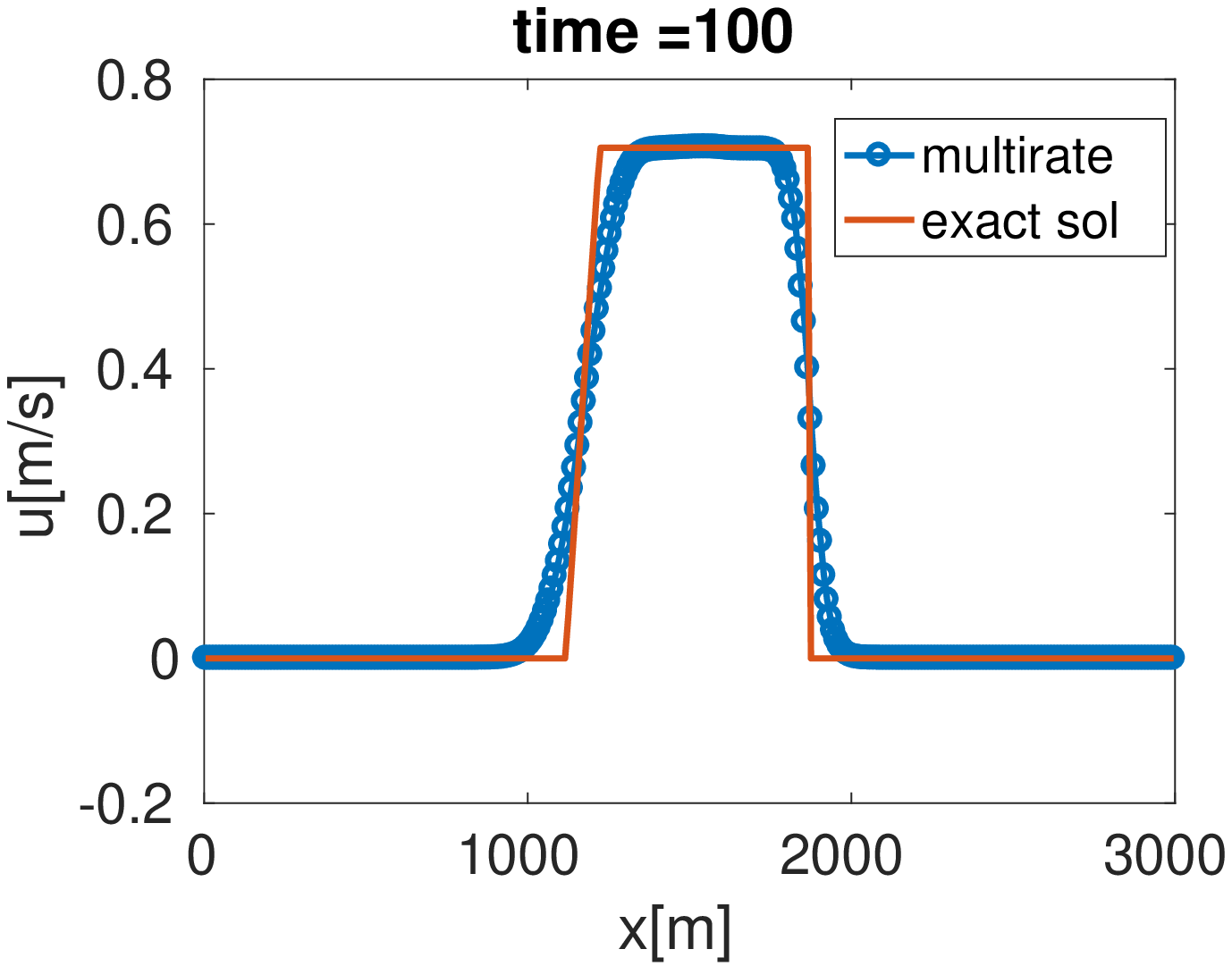}
\caption{Solutions at time $t=42$ and $t = 100$ for the $h$ variable (on the left) and for $u = \frac{q}{h}$ variable (on the right).}
\label{dam:sol}
\end{figure}

\begin{figure}
\centering
\includegraphics[width=0.6\textwidth]{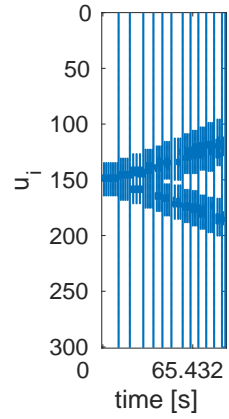}
\caption{Set of active components for one variable at each time.}
\label{dam:mesh}
\end{figure}

\begin{figure}
\centering
\includegraphics[width=0.6\textwidth]{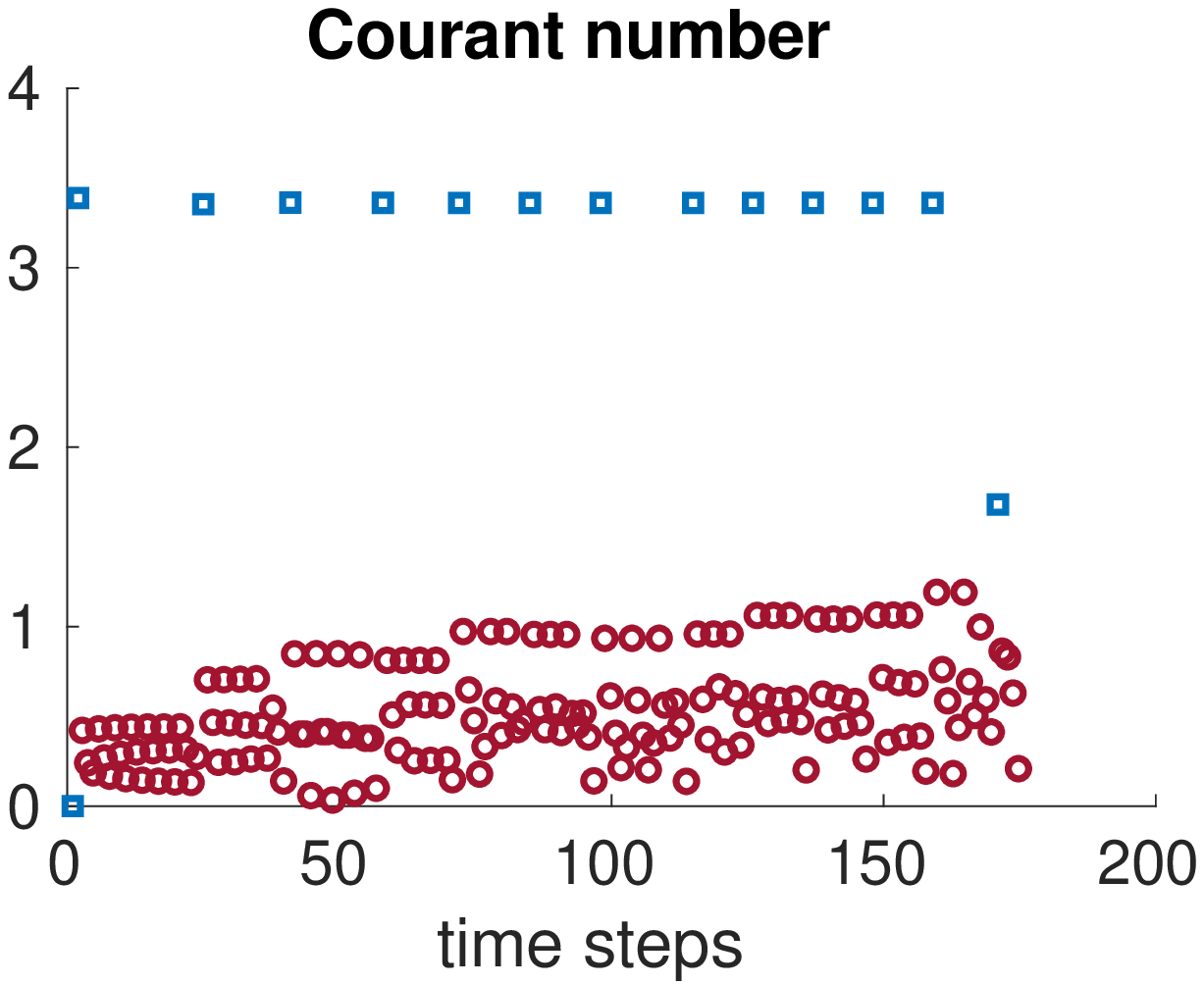}
\caption{Courant number for each time step for the dam break problem.}
\label{dam:cour}
\end{figure}

When performing  this test with the original version of the algorithm described
in the previous sections,  numerical oscillation across the boundary between the refinement and the non-refinement regions were observed. These oscillations are due to the fact that the error estimator accepted some fluxes that were changing their values inside the time slab and it was not correct to use their final time slab values for the entire considered sub-step.
To avoid this problem, we slightly modified the set of rejected fluxes. If a flux is rejected, we also reject a number of fluxes (on the left or on the right or on both sides, depending on the sign of the eigenvalues) equal to the local Courant number. In this way, as shown in Fig. \ref{dam:sol}, the solution has the correct behavior; of course, we are increasing the set of active components, but the latent components are still the majority during the time integration (Fig. \ref{dam:mesh}).  
It can be seen clearly that, as in the scalar case, the method is able to identify automatically the complex nonlinear features of the flow.    It can also be seen in Fig. \ref{dam:cour} that a Courant number larger than one was allowed
for the global time steps without any significant loss in accuracy.

\subsection{Shallow water equations with rotation}
We have also considered the shallow water equations with rotation, which are a classical
idealized model for the phenomenon of geostrophic adjustment, see e.g.~\cite{gill:1982}.
This system, in the semi-linear form obtained discarding the nonlinear momentum advection terms,
can be written as:

\begin{equation}
\begin{cases}
\frac{\partial \eta}{ \partial t} + \frac{\partial ((\eta + \eta_0)u)}{\partial x} = 0  & (x,t) \in(-L,L)\times (0,T)\\
\frac{\partial u}{\partial t} + g\frac{\partial \eta}{\partial x} + f v = 0 & (x,t)\in(-L,L)\times (0,T)\\
\frac{\partial v}{\partial t} -f u=0 & (x,t) \in(-L,L)\times (0,T)\\
\eta(x,t=0) = \exp\left(-\frac{(50x)^2}{(2L)^2}\right) &x\in(-L,L)\\
u(x,t=0) = v(x,t=0) = 0 &x \in(-L,L)\\
\eta(-L,t)= \eta(L,t)=0 & t \in (0,T)\\
u(-L,t)= u(L,t)=0 & t \in (0,T)\\
v(-L,t)= v(L,t)=0 & t \in (0,T)\\
\end{cases}
\end{equation}
Here, $\eta $ denotes the free surface height, $u$ the velocity in the $x $ direction, 
$g $ is the acceleration of gravity, $f$ is a constant Coriolis parameter
and  $v$  represents the velocity in the direction orthogonal to the one dimensional flow being considered. 
This system is of particular interest since it describes a dynamics with two different time scales,
a fast one associated to the propagation of external gravity waves and a slow one associated
with rotational effects  and the onset of geostrophic equilibrium. Semi-implicit techniques commonly applied
for geophysical scale flows (see e.g. the classical paper \cite{robert:1982} and
\cite{giraldo:2013}, \cite{tumolo:2015} for two more modern examples of this approach)
allow to achieve an accurate approximation of the slow components, while sacrificing the accuracy 
of the fast ones.

In order to represent a large geophysical scale, we have used $L=8\times10^{6}$ m, $T = 3\times10^{6}$ s, $f = 1\cdot10^{-4}$ $1/$s and $\eta_0 = 1000$ m.
We have discretized in space with $N_x = 480$ cells and we have used, as space discretization, the conservative   centered finite difference scheme:

\begin{align*} 
& \frac{d\eta_i}{dt} = - \left[ \frac{u_i\eta_i + u_{i+1}\eta_{i+1}}{2\Delta x}  - \frac{u_i\eta_i+u_{i-1}\eta_{i-1}}{2\Delta x}\right],\\   
&\frac{du_i}{dt} = - g \left[ \frac{\eta_i + \eta_{i+1}}{2 \Delta x}- \frac{\eta_i + \eta_{i+1}}{2\Delta x} \right] - fv_i,\\
&\frac{dv_i}{dt} = fu_i.
\end{align*}  

\begin{figure}
\centering
\includegraphics[width=0.6\textwidth]{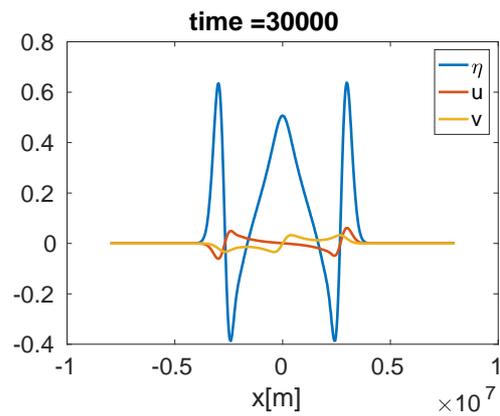}
\caption{Solutions at the final time computed with the multirate method.}
\label{sw:sol}
\end{figure}

\begin{figure}
\centering
\includegraphics[width=0.84\textwidth]{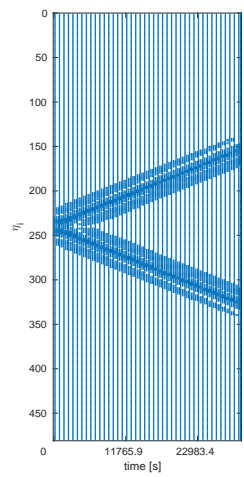}
\caption{ Set of active components for $\eta$ variable.}
\label{sw:mesh}
\end{figure}

\begin{figure}
\centering
\includegraphics[width=0.6\textwidth]{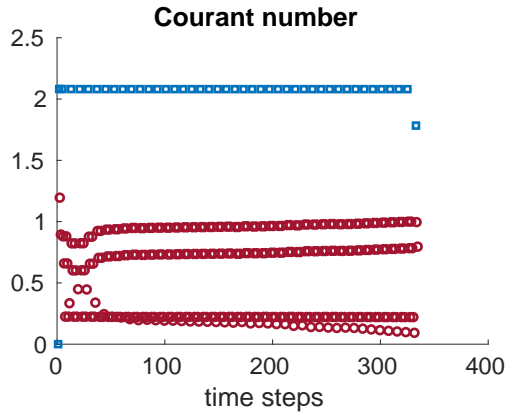}
\caption{Courant number for each time step for the shallow water equation with rotation.}
\label{sw:courant}
\end{figure}

In this case,  we used a global step $\Delta t= 700 s$  to discretize in time. 
The solution is represented in Fig.~\ref{sw:sol}, while
the set of active/refined components for the  $\eta$ variable is displayed in Fig. \ref{sw:mesh}. 
It can be seen that, also in this case, the proposed algorithm is able to identify automatically
the different time scales present in the solution. The component of the solution at the center of the domain,
which tends to geostrophic equilibrium on a slow time scale, does not require any refinement of the time step,
while the fast propagating gravity waves induce refinement along the wave trails.
Notice that we plot the active components   for the $\eta$ variable only because, as explained in section 
\ref{system_ref}, the set of active components and active fluxes are the same for each variable of the system in order to preserve mass.  It can also be seen in Fig. \ref{sw:courant}   that Courant numbers larger than one are feasible for the global time steps without any significant loss in accuracy.  
   
\begin{table}
\caption{Computational time, number of time steps and total components number involved, using a relative tolerance equal to $ 1e-4$ and as absolute tolerance $1e-3$ for both the single rate and the multirate approach.}
\centering
\label{sw-table}
\begin{tabular}{llll}
\hline\noalign{\smallskip}
 & comp. time [s]  & $\#$ time steps & $\#$ function eval. \\
\noalign{\smallskip}\hline\noalign{\smallskip}
Multirate		  & $74.71$ & $424 (83 \mbox{ glob steps})$ &$102336$\\
\noalign{\smallskip}\hline\noalign{\smallskip}
Single rate & $179.29$ & $173$& $186810$ \\
\noalign{\smallskip}\hline
\end{tabular}
\end{table}

In Table \ref{sw-table} we reported the comparison with the single-rate version of the TR-BDF2 method. In the first column we report  the CPU time required to solve the problem with the two different methods. In the second column we report the number of time steps necessary with each approach until final time. It is to be remarked that
both methods were implemented in a rather straightforward way and that the respective codes are far from optimized. On the other hand, exactly the same computational components, such as e.g. the Newton solver,
were employed in both, so that the ratio of the CPU times required by the two approaches is a reasonable estimate
of the potential speed-up. It can be seen that  the multirate approach solves the problem more than twice as fast than the single rate method.

The multirate method uses more time steps with respect to the the single rate method, but only  roughly $20\%$ of these are global time steps,
while for the remaining   time steps   only few components have to be computed. In fact, in the third column of the table we report the number of components involved to solve the system from the initial time to the final time. The single-rate method involves about twice as many components as the multirate method.

\section{Conclusions}
We propose a conservative implicit multirate method for time integration of hyperbolic problems. To integrate in time we have used the TR-BDF2 method, which is a second order, L-stable implicit method, but the approach can be easily generalized to other implicit methods. 

The partition of fast and slow components is based on the numerical flux, in order to preserve the conservative
nature of the spatial discretizations employed. A consistency
 analysis  has been carried out, showing that only implicit discretizations  that do not involve previous values
 in the computation of the fluxes, such as the backward Euler method, are fully consistent. On the other hand,
 inconsistency only arises at the interface between refined and non refined regions and does not seem
 to affect the accuracy of the method significantly.

We have tested this approach on several scalar equations and, to the best of out knowledge for the first time,
we have applied a self-adjusting multirate method to systems of non-linear conservation laws, albeit
 only in the one dimensional case.  
The results show that the multirate approach captures automatically the behaviour of the solution and refines only where it is necessary, thus achieving a reduction of the CPU costs without significant losses of accuracy. 
The extension of this method to more complex problems and to multi-dimensional equations is an area of current research.       

\begin{acknowledgements}
The first, second and fourth authors would like to acknowledge the financial support of the 
INDAM - GNCS  projects {\it Metodi numerici semi-impliciti e semi-Lagrangiani  per sistemi iperbolici di leggi di bilancio (2015)} (second author only) and {\it Modellazione numerica di fenomeni idro /geomeccanici per la simulazione di eventi sismici (2017)}.
\end{acknowledgements}

\bibliographystyle{spbasic}      


\end{document}